\documentclass[journal abbreviation, manuscript]{copernicus}

\begin{document}
\nolinenumbers

\title{In-Flight Estimation of Instrument Spectral Response Functions Using Sparse Representations}

\Author[1,2]{Jihanne}{El Haouari}
\Author[3]{Jean-Michel}{Gaucel}
\Author[4]{Christelle}{Pittet}
\Author[1,2]{Jean-Yves}{Tourneret}
\Author[2]{Herwig}{Wendt}

\affil[1]{TéSA laboratory, Toulouse, France}
\affil[2]{Univ. Toulouse, IRIT-ENSEEIHT, Toulouse, France}
\affil[3]{Thales Alenia Space Cannes, France}
\affil[4]{Centre National d'Etudes Spatiales, Centre Spatial de Toulouse, France}

\correspondence{Jihanne El Haouari (jihanne.elhaouari@prd.tesa.fr)}

\runningtitle{In-Flight Estimation of Instrument Spectral Response Functions Using Sparse Representations}

\runningauthor{J. El Haouari, J.-M. Gaucel, C. Pittet, J.-Y. Tourneret, H. Wendt}

\received{}
\pubdiscuss{} 
\revised{}
\accepted{}
\published{}

\firstpage{1}

\maketitle

\begin{abstract}
Accurate estimates of Instrument Spectral Response Functions (ISRFs) are crucial in order to have a good characterization of high resolution spectrometers. Spectrometers are composed of different optical elements that can induce errors in the measurements and therefore need to be modeled as accurately as possible. Parametric models are currently used to estimate these response functions. However, these models cannot always take into account the diversity of ISRF shapes that are encountered in practical applications. This paper studies a new ISRF estimation method based on a sparse representation of atoms belonging to a dictionary. This method is applied to different high-resolution spectrometers in order to assess its reproducibility for multiple remote sensing missions. The proposed method is shown to be very competitive when compared to the more commonly used parametric models, and yields normalized ISRF estimation errors less than $1\%$. 
\end{abstract}

\introduction  
Space remote sensing makes it possible to remotely measure the composition of the atmosphere or the troposphere and to retrieve trace gases concentrations. It can also be used to monitor molecule fluxes at the Earth's surface, such as for the MicroCarb mission that is designed to monitor $CO_2$ fluxes \citep{Cansot2022}. This can be done by analyzing the interaction of the atmosphere with natural radiation, such as the sunlight, or artificial radiation, generated for example by a laser. Indeed, the presence of some molecules in the path of radiation modifies its spectral content at the characteristic wavelengths of the different elements. The information directly obtained from satellites is the atmospheric spectrum. By considering some specific wavelengths of interest, it is possible to determine the concentration of the desired trace gases in a column of atmosphere by comparing these measured spectra with a reference spectrum obtained using a radiative transfer model. 

The instruments used for gas concentration estimation are high resolution spectrometers. This study focuses on six different spectrometers: Avantes \citep{Avantes} that is used for fluorescence analysis or reflexion spectrometry in the UV/NIR ranges, GOME-2 \citep{Munro2016} and OMI \citep{Dobber2006} that mainly monitor atmospheric ozone at UV-visible wavelength ranges, TROPOMI \citep{Kleipool2018} that measures the UV, visible, near-infrared and shortwave infrared reflectance of the Earth, OCO-2 \citep{Lee2017} and MicroCarb \citep{Cansot2022} that are dedicated to study the atmospheric carbon dioxyde and oxygen in order to determine their concentrations at the Earth surface. The analysis of these concentrations for the two last spectrometers can provide a better understanding of the carbon cycle, which is important in the context of climate change.

Spectrometers consist mainly of an optical part (for example composed of a slit, a telescope and dispersive grating) and a detector. In this configuration, the telescope projects the image of the Earth onto the spectrometer slit and then onto the detector. Each pixel of the detector is associated with a spatial direction (called ACT for ACross Track) and a specific wavelength. A binning and an averaging  along the ACT axis are performed in order to improve the Signal to Noise Ratio (SNR). For each of the two parts (optical part and detector), a response function is defined, which leads to a continuous optical function and another function associated with each pixel of the detector. This results in a global response function associated with each pixel along the spectral axis, known as the Instrument Spectral Response Function (ISRF), associated with a specific wavelength. The ISRFs can vary significantly depending on the instrument considered and their shapes depend on the central wavelength, among other factors. The estimation of trace gas concentrations is an inversion process that is performed on the ground from spectrometer measurements and the instrument ISRFs. The accuracy of this estimation highly depends on the knowledge of these ISRFs for all pixels. For some missions, ISRFs are expected to be known with a normalized error less than $1\%$, which is a challenge since wavelength variations of ISRFs exceed this threshold in most cases.

Spectrometers are first calibrated on ground where their associated ISRFs are estimated experimentally. However, the ISRFs are subject to in-flight changes due to mechanical movements associated with the launch of the instruments, thermal changes in orbit, or certain sensitivities linked to the instrument itself (such as MicroCarb's sensitivity to the scene). As a consequence, these ISRFs need to be re-estimated regularly in practical applications in-flight throughout the mission. This implies a real time process and therefore it is important to have an accurate and fast enough ISRF estimation method.
 The principle of the estimation is to take a measurement of a spectrally known scene and to compare it with a spectral model of the scene convolved with the ISRFs at different wavelengths. 
Parametric models have been widely used in the literature to estimate ISRFs. Gaussian and generalized Gaussian parametric models (referred to as ``Gauss'' and ``Super-Gauss'') were proposed in \cite{Beirle2017}.Parametric models are attractive for their simplicity and small number of parameters. However, they are not flexible enough to represent the range of ISRF shapes adequately. The ISRF estimation problem and the most important parametric models that have been considered in the literature are detailed in Section \ref{section:ISRF_models_and_literature_methods}.

The objective of this work is to overcome the limitations of the existing parametric ISRF estimation methods caused by their insufficient accuracy. To that end, we propose as a first major contribution a new estimation strategy based on sparse representations of the ISRFs in a dictionary of well chosen atoms. More precisely, the ISRFs are decomposed in a dictionary that is constructed using several ISRFs that are available from ground characterization for each instrument. The dictionaries can also be updated iteratively on-line. For each instrument, each ISRF is then approximated by a linear combination of a small number of atoms of the dictionary associated with the instrument. The proposed approach will be referred to as SPIRIT, for ``SParse representation of Instrument spectral Response functions using a dIcTionary'', and is detailed in Section \ref{section:Sparse_approximation}. We investigate and compare two different methods for obtaining the sparse representations of ISRFs. 
As a second contribution, we conduct an extensive numerical study of the proposed ISRF estimation approach and compare it to parametric methods for datasets from six different spectrometers used in space missions, whose characteristics are detailed in Section \ref{section:Instruments_and_datasets}. Numerical results for these datasets, as well as results for estimation performance and for robustness with respect to design choices and noise corrupting the observed measurements are reported in Section \ref{section:Results_and_discussion}. 
As a major outcome, the proposed SPIRIT method yields significantly improved flexibility and accuracy for ISRF estimation when compared to previous state-of-the-art parametric methods, consistently through the different datasets and scenarios, with a small number of parameters that can be easily and efficiently estimated in real-time.

\section{Existing models and estimation methods}
\label{section:ISRF_models_and_literature_methods}

\subsection{ISRF estimation model}
 The ISRF, that is sometimes referred to as Instrument Line Shape (ILS) \citep{Sun2017b} or Slit Function \citep{Sun2017a}, is a function that describes the response of an instrument to a given wavelength. In this work, we only consider the spectral information and thus each ``pixel'' $l$ is associated with a specific wavelength $\lambda_l$ yielding an ISRF at this wavelength.
The in-flight identification of ISRFs is obtained from scenes that are assumed to be perfectly known  radiometrically and spectrally (such as Sun, Moon, uniform scenes such as desert, etc.), which are referred to as reference spectra. The principle of ISRF estimation is to determine the in-flight ISRFs for each wavelength $\lambda_l$ that minimize some similarity measure between the measured spectrum  (denoted as $s(\lambda_l)$) and the reference spectrum (denoted as $r(\lambda_l)$) convolved with the ISRF (denoted as $I_l(\lambda_l)$:
 \begin{align}
 s(\lambda_l) = (r * I_l)(\lambda_l) = \int_{\lambda_{\text{min}}}^{\lambda_{\text{max}}} r(\lambda_l - u) I_l(u) du, \quad l=1,...,N_{\lambda}.
\end{align}  
$*$ denotes convolution and $N_{\lambda}$ is the number of central wavelengths $\lambda_l$, each associated with one ISRF $I_l$.

For practical purposes, this equation can be discretized leading to:
   \begin{align}
   \label{sect2.1:Optimization_model}
 s(\lambda_l) \approx \sum_{n = - N/2 } ^{N/2}  r(\lambda_l - n \Delta) I_l(n \Delta), \quad l=1,...,N_{\lambda},
\end{align}  
where $\Delta$ is the wavelength sampling interval for the ISRFs, which is assumed to be regularly sampled. In other words, a vector $\boldsymbol{I}_l = [I_l(-\frac{N}{2} \Delta) ,... ,   I_l(\frac{N}{2} \Delta)]^T \in \mathbb{R}^{N + 1}$ needs to be estimated for each ISRF, corresponding to the values that it takes on the wavelength grid $\boldsymbol{\Delta} =[ - \frac{N}{2} \Delta, ... ,\frac{N}{2} \Delta ]^T \in \mathbb{R}^{N + 1}$.
  The objective of the ISRF estimation problem is to solve the inverse problem \eqref{sect2.1:Optimization_model} assuming knowledge of both the reference spectrum $r$ and the measurements $ s(\lambda_l)$.
 
A major difficulty with the inverse problem \eqref{sect2.1:Optimization_model} is that there is only one measurement per fixed wavelength, which makes it impossible to estimate the vector $\boldsymbol{I}_l$ without further assumptions. Two approaches can be used to make this estimation problem identifiable. The first idea is to consider knowledge of several reference spectra $r$ for every wavelength. The problem is that this would not only require a sufficient number of calibration scenes to be available, but also that they substantially differ for each wavelength in order to provide complementary information on the shapes of the ISRFs.
The second method, which is considered in this paper, makes use of only one reference spectrum and is based on the assumption that the ISRFs for adjacent wavelengths $\lambda_l$ are similar. This is a reasonable assumption for the ISRFs of real-world spectrometers.  Then, to estimate an ISRF at wavelength $\lambda_l$, we propose to consider a vector $\boldsymbol{s}_l = [s(\lambda_{l - \frac{N_{\text{obs}}}{2}}), ... , s(\lambda_{l + \frac{N_{\text{obs}}}{2}}) ]^T \in \mathbb{R}^{N_{\text{obs}} + 1} $ of $N_{\text{obs}}+1$ observations, including also those from the neighboring ISRFs.
Rewritten in matrix form, \eqref{sect2.1:Optimization_model} simplifies to:
$$
\boldsymbol{s}_l = \boldsymbol{R}_{l} \boldsymbol{I}_l ,
$$
where $\boldsymbol{R}_{l} = [\boldsymbol{r}_{ l - \frac{N_{\text{obs}}}{2}}, ... ,\boldsymbol{r}_{l + \frac{N_{\text{obs}}}{2}} ]^T \in \mathbb{R}^{(N_{\text{obs}} + 1) \times (N + 1)} $ contains the values $\boldsymbol{r}_{l} = [r(\lambda_{l} - \frac{N}{2} \Delta), ... , r(\lambda_{l} + \frac{N}{2} \Delta)] \in   \mathbb{R}^{N + 1}$ of the reference spectrum covered by the different ISRFs in the neighborhood.
Given a model for the ISRF, estimating $\boldsymbol{I}_l$ can then be conducted for each wavelength $\lambda_l$ by minimizing  the residual error $|| \boldsymbol{s}_{l} - \boldsymbol{R}_{l} \boldsymbol{I}_l ||_2^2$.
 
\subsection{Parametric models}

A classical way to model and estimate the ISRF at wavelength $\lambda_l$ is to use a parametric Gaussian model defined by:
\begin{align} \label{Gauss}
	\boldsymbol{I}_{l,\boldsymbol{\beta}_{\text{G}} }(x) =  A_G \exp\left[- \frac{( \lambda_l - x - \mu_{\text{G}})^2  }{2\sigma_G^2} \right], \quad l=1,...,N_{\lambda},\quad x\in \boldsymbol{\Delta},
\end{align}
where $\boldsymbol{\beta}_G = [A_{\text{G}}, \mu_{\text{G}},\sigma^2_G]^T $ is the unknown vector of parameters to be estimated.

\citet{Beirle2017} studied an alternative ISRF model using a generalized Gaussian distribution referred to as ``super-Gaussian'' in order to better fit the ISRF shapes: 
\begin{align} \label{SuperGauss}
	\boldsymbol{I}_{l,\boldsymbol{\beta}_{\text{SG}} }(x) =  A_{\text{SG}} \exp\left[- \left| \frac{ \lambda_l - x - \mu_{\text{SG}}  }{w_{\text{SG}}} \right|^{k_{\text{SG}}} \right], \quad l=1,...,N_{\lambda},,\quad x\in \boldsymbol{\Delta},
\end{align}
where $\boldsymbol{\beta}_{\text{SG}} = [A_{\text{SG}}, \mu_{\text{SG}},w_{\text{SG}},k_{\text{SG}}]^T $ is the  unknown parameter vector to estimate. This model reduces to the Gaussian model when $w_{\text{SG}} = 2 \sigma^2_{\text{G}}$ and $k_{\text{SG}} = 2$. The parameters $w_{\text{SG}}$ and $k_{\text{SG}}$ are the scale and shape parameters of the distribution, allowing more or less flattened ISRFs to be obtained.

When using the parametric models \eqref{Gauss} and \eqref{SuperGauss}, the ISRF estimation problem consists of estimating the unknown model parameters for each sliding window. This estimation can be performed using the least squares method, which minimizes the following cost function:
\begin{equation}
	C_l(\boldsymbol{\beta}) = \sum_{n=1}^{N+1} = ||\boldsymbol{s}_l - \boldsymbol{R}_{l} \boldsymbol{I}_{l, \boldsymbol{\beta}} ||_2 ^2, \quad l=1,...,N_{\lambda},
\end{equation}
where $\boldsymbol{\beta} \in \{ \boldsymbol{\beta}_{\text{G}} , \boldsymbol{\beta}_{\text{SG}} \}$ is the unknown parameter vector and $\boldsymbol{I}_{l, \boldsymbol{\beta}} =[\boldsymbol{I}_{l,\boldsymbol{\beta}} (\delta_1),...,\boldsymbol{I}_{l,\boldsymbol{\beta}} (\delta_{N+1})]^T $. 

\begin{figure}[h!] 
            \includegraphics[trim={0cm 0cm 0cm 0cm}, clip,width=10.0cm]{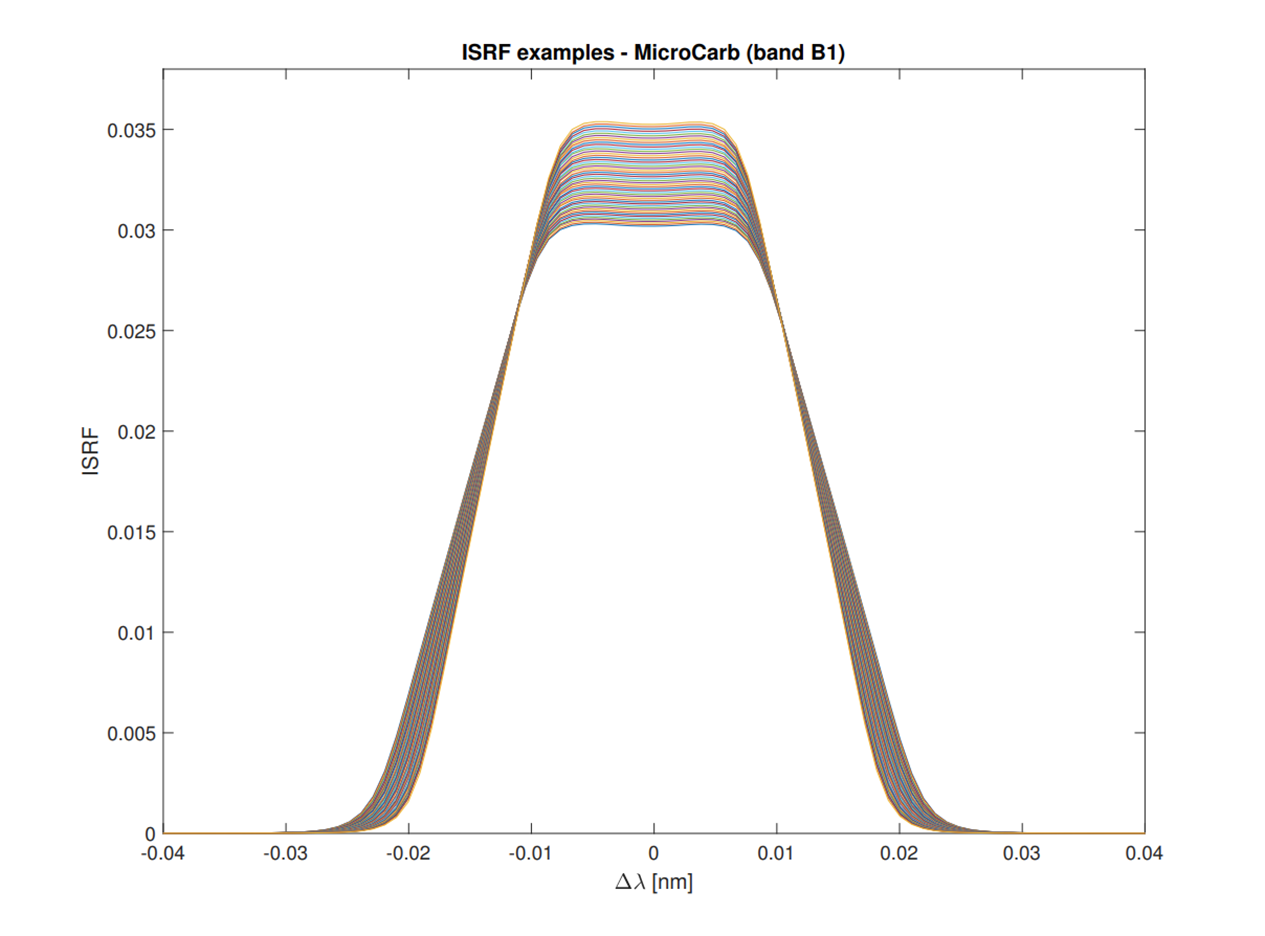}
    \caption{Examples of MicroCarb ISRFs.}
    \label{figure:ISRF_example_MC}
\end{figure}

It is difficult to analytically construct accurate forward models with a small number of parameters for ISRFs because they would need to incorporate a significant number of ``contributors'' associated with the instrument optics (slit, mirror, lens, separator, dispersing element), the detector or the acquisition mode.
The state-of-the-art therefore considers simple parametric models, such as Gaussian or generalized Gaussian models, that can struggle to take into account the variety of different ISRF shapes that can be observed in practice. 
An illustration is provided in Fig. \ref{figure:ISRF_example_MC}, which shows examples of ISRFs for the MicroCarb mission. Clearly, these ISRFs cannot be accurately modeled by bell-shaped Gaussian distributions and by generalized Gaussians (e.g., because of the dip at the center). Other examples for different instruments are displayed in Fig.~\ref{figure:illustration_single_ISRF}. This motivates the study of  a new estimation method for ISRFs. 
\section{Sparse approximations of the ISRFs}
\label{section:Sparse_approximation}

In this paper, we investigate the use of sparse representations for ISRFs in a dictionary of well chosen atoms. Models based on sparse approximations and on dictionary learning have been widely and successfully used for different signal and image processing applications \citep{SurveySparseRep}. These applications include image denoising, image classification, image reconstruction, compressed sensing or dimensionality reduction and can model large varieties of signals and images \citep{Gradient_Projection,dictionary_learning}. However, their use has never been reported before in the context of spectrometers, ISRF and instrument parameter estimation, which is precisely the objective of this work.
 
\subsection{Construction of the dictionary}
\begin{figure}[h!] 
        
            \includegraphics[trim={4cm 0.5cm 4cm 0.5cm}, clip,width=17.0cm]{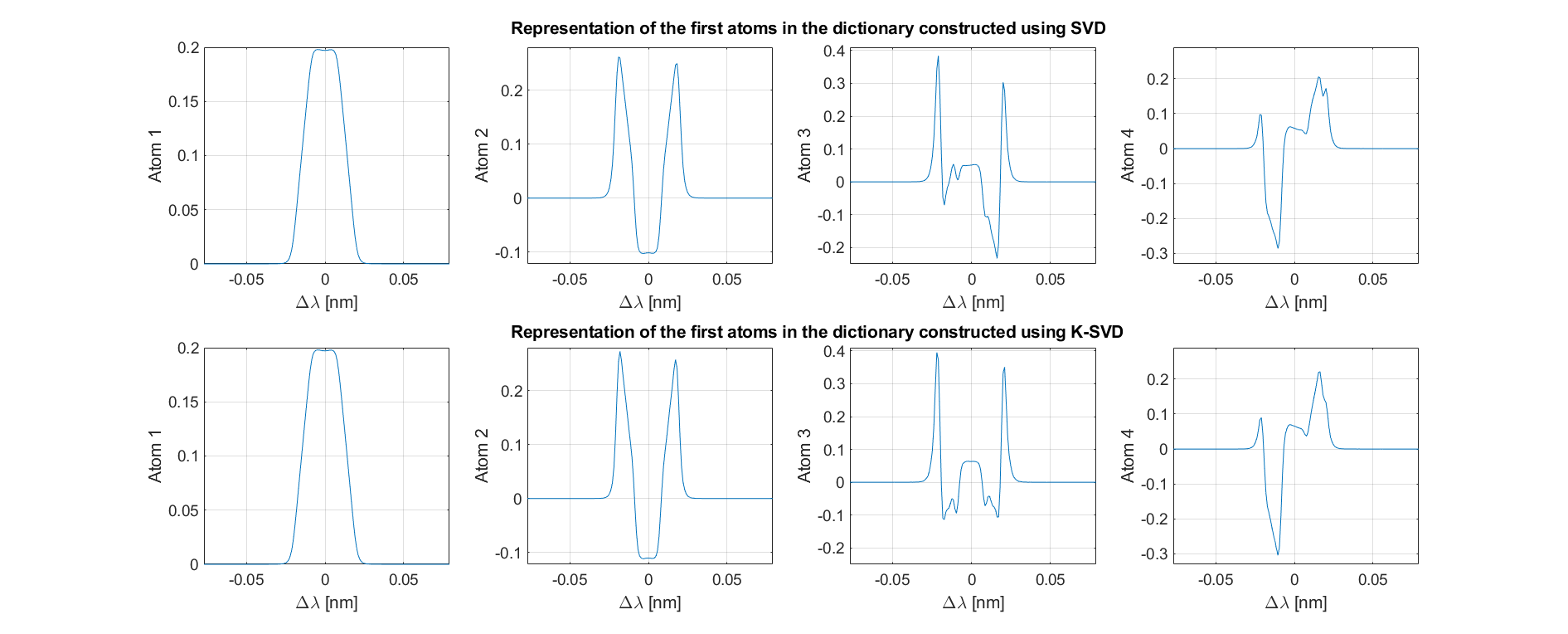}
    \caption{Representation of the four first atoms of the dictionary constructed using one SVD (top) or using the K-SVD algorithm (bottom) for the MicroCarb spectrometer (band B1).}
    \label{figure:dictionary_atoms}
\end{figure}

Sparse representations consist in expressing a given signal as a linear combination of a small number of signals that belong to a collection of reference patterns, or atoms, which is called a dictionary. This paper proposes to decompose the ISRF in a dictionary of atoms $\boldsymbol{\Phi} \in \mathbb{R}^{(N + 1) \times N_{\text{D}}}$: 
\begin{align} \label{ISRF}
	\boldsymbol{I}_l \approx \boldsymbol{I}_l^K = \boldsymbol{\Phi} \boldsymbol{\alpha}_l = \sum_{k=1}^K \boldsymbol{\Phi}_{\gamma_k} \alpha_{l,k}, \quad l=1,...,N_{\lambda},
\end{align}
where $ \boldsymbol{\Phi}_{\gamma_k}$ is the $\gamma_k$th selected atom, i.e., the $\gamma_k$th column of the dictionary $\boldsymbol{\Phi}$ and $\alpha_{l,k}$ is the associated non zero coefficient of the sparse vector $ \boldsymbol{\alpha}_l=[\alpha_{l,1},...,\alpha_{l,K}]^T \in \mathbb{R}^{N_{\text{D}}}$. The dictionary is built in such a way that its atoms (i.e., its columns) provide an efficient representation of the signal. Different methods allowing the dictionary to be built have been proposed in the literature. These methods are based on probabilistic learning, clustering, vector quantization or Bayesian inference \citep{dictionary_learning}. Dictionary learning usually involves a two-stage optimization structure, consisting first of a sparse coding step, to find the sparse vector $\boldsymbol{\alpha}_l$ which minimizes the objective function $||\boldsymbol{I}_l - \boldsymbol{\Phi} \boldsymbol{\alpha}_l ||_2^2 $ for a fixed dictionary $\boldsymbol{\Phi}$ and then a dictionary update step, where the dictionary is estimated for a fixed sparse vector $\boldsymbol{\alpha}_l$. Depending on the application, the dictionary can be updated using a closed form solution, gradient descent, or using ground truth data. In this work we investigate two different ways of building the dictionary $\boldsymbol{\Phi}$. The first method constructs $\boldsymbol{\Phi}$ by using the $N_{\text{D}}$ singular vectors associated with the largest singular values of the SVD of a matrix composed of representatives ISRF examples. The second method uses the K-SVD algorithm of \citet{Elad2006}, which belongs to the state-of-the-art. The K-SVD algorithm is a generalization of the K-means algorithm in which the dictionary is updated by changing its columns separately and sequentially and applying $K$ singular value decompositions (SVDs) on an appropriate error matrix.  Fig.~\ref{figure:dictionary_atoms} displays the first atoms of dictionaries constructed using these two methods for the band B1 of MicroCarb. These dictionaries are similar, especially the two first atoms that correspond to the most energetic singular values. The two first atoms can be interpreted as the approximate average of all ISRFs used to build the dictionary (first atom), and a correction for adjusting the different widths of the ISRFs for different wavelengths (second atom), as seen in Fig.~\ref{figure:ISRF_example_MC}. The higher order atoms obtained with SVD and K-SVD are slightly different but with the same global shape.

\subsection{Inverse problem}
Assuming that the ISRF can decomposed in the dictionary $\boldsymbol{\Phi}$ as in \eqref{ISRF}, the measured spectrum can be written as follows: 
\begin{align*}
	\boldsymbol{s}_{l} \approx \boldsymbol{R}_{l} \boldsymbol{I}_l \approx \boldsymbol{R}_{l} \boldsymbol{\Phi} \boldsymbol{\alpha}_l =\boldsymbol{\Psi}_l \boldsymbol{\alpha}_l, \quad l=1,...,N_{\lambda}.
\end{align*}
Thus, the ISRF estimation problem reduces to finding the sparse vector $\boldsymbol{\alpha}_l$ that minimizes the residual $||\boldsymbol{s}_{l} - \boldsymbol{\Psi}_l \boldsymbol{\alpha}_l ||_2^2$. 
This sparse coding problem has been mathematically formulated in different ways \citep{SurveySparseRep}. One can use the $l_0$ pseudo-norm regularization$ ||.||_0$ with a penalty parameter $\mu$, leading to the following problem:
\begin{align}
\label{Optimization_problem}
	\arg \min_{\boldsymbol{\alpha}_l} L(\boldsymbol{\alpha}_l, \mu) = \arg \min_{\boldsymbol{\alpha}_l} ||\boldsymbol{s}_{l} - \boldsymbol{\Psi}_l \boldsymbol{\alpha}_l ||_2^2 + \mu ||\boldsymbol{\alpha}_l||_0, \quad l=1,...,N_{\lambda}.
\end{align}
This problem being non-convex and NP-hard, many approximations and heuristics have been proposed in the literature to find an approximate solution. 
A standard method consists of using greedy algorithms such as the Orthogonal Matching Pursuit (OMP). OMP is a modification of the Matching Pursuit (MP) algorithm, which improves convergence by adding an orthogonalization step \citep{Mallat1993,Pati1993}. The atoms of the dictionary best approximating the data fidelity term $||\boldsymbol{s}_{l} - \boldsymbol{\Psi}_l \boldsymbol{\alpha}_l ||_2^2$, are iteratively determined by minimizing the remaining residual error.
Another method replaces the pseudo-norm $l_0$ in \eqref{Optimization_problem} by the $l_1$ norm, which leads to a convex problem known as the LASSO problem \citep{LASSO_PART1}:
\begin{align}
\label{lasso}
	\arg \min_{\boldsymbol{\alpha}_l} L(\boldsymbol{\alpha}_l, \mu) = \arg \min_{\boldsymbol{\alpha}_l} ||\boldsymbol{s}_{l} - \boldsymbol{\Psi}_l \boldsymbol{\alpha}_l ||_2^2 + \gamma ||\boldsymbol{\alpha}_l||_1, \quad l=1,...,N_{\lambda},
\end{align}
and the related algorithms studied in, e.g., \citep{Gradient_Projection,Interior-Point}. 

These algorithms provide a highly flexible decomposition of the ISRF, as the choice of the dictionary is not constrained to a specific form. 
Indeed, the basis functions that can be learned, for example by using the K-SVD algorithm in conjunction with various Matching Pursuit algorithms. 
This paper compares the use of fixed dictionaries obtained by a single SVD and dictionaries estimated by alternating between an SVD step to update the dictionary and an OMP step for sparse coding (K-SVD). The proposed approach using OMP or Lasso (or other sparse formulations) and either fixed or re-estimated dictionaries will be referred to as SPIRIT for ``SParse representation of Instrument spectral Response Functions using a dIcTionary''.

\section{Instruments, datasets \& preprocessing}
\label{section:Instruments_and_datasets}

This section considers different high-resolution spectrometers in order to assess the robustness of the proposed ISRF estimation method to different instruments used in space missions. 

\subsection{Synthetic data generation}
Reference spectra used in this study were generated using the 4A/OP (Automatized Atmospheric Absorption Atlas) software described in \citet*{4AOP}. This software is based on a fast and accurate line-by-line transfer model that can be integrated in operational processing chains including inverse problem processing \citep{4AOPPresentation}. It was selected as the official radiative model and reference code by CNES for the MicroCarb mission \citep{4AOPPresentation}. Typical ISRFs were recovered and then normalized to area 1 for each instrument (see details in the next paragraphs), convolved to the reference spectra and embedded in additive Gaussian noise to generate representative measurements. The advantage of this data generation method is to provide ground truth ISRFs, which can be used to assess the performance of the different methods in a controlled scenario. 

\subsection{Avantes}
Avantes \citep{Avantes} is an ultra-low straylight fiber spectrometer composed of standard $2048 \times 64$ back-thinned CCD detectors. This spectrometer has multiple uses and can cover different ranges in the UV, VIS or NIR from 200 to 1160 nm. The spectrometer resolution varies between $0.09$ nm and $20$ nm depending on its configuration. Avantes has a single detector of $2048$ spectral pixels and covers the wavelength range 296-459 nm. The data used for this spectrometer was measured in April 2015 and $7$ ISRFs were made available in the supplement material of \citet{Beirle2017}. These $7$ templates were interpolated using splines to generate $N_{\lambda}=1741$ frequency varying  ISRFs with $N = 199$ samples.

\subsection{GOME-2}
The Global Ozone Monitoring Experiment-2 (GOME-2) is an optical spectrometer that was launched in the Metop-A EUMETSAT satellite in october 2006. All the information on the design and characterization of the spectrometer are described in \citet{Munro2016}. The aim of the GOME-2 mission is to obtain information on trace gas concentration by sensing the Earth's backscattered radiance in the UV part of the spectrum (240-790 nm). The instrument is composed of four main optical channels. For each channel there is one detector of $1024$ pixels, with a high resolution ranging from 0.26 nm to 0.51 nm. This article concentrates on the third band of GOME-2 (397 nm-604 nm). Data for the ISRFs were made available by EUMETSAT (ftp://ftp.eumetsat.int/pub/EPS/out/GOME/Calibration-Data-Sets/Slit-Function-Key-Data/FM3-Metop-A/). The ISRFs were obtained in January 2007 using the FM3 model \citep{Siddans2018} and were interpolated using splines to generate ISRFs with $N_{\lambda}=877$ frequencies and a sample size $N =705$.

\subsection{OMI}
The Ozone Monitoring Instrument (OMI) is an hyperspectral imager onboard the NASA's Earth Observing system satellite Aura launched in July 2004. The aim of the mission is to determine trace gases at high spectral and spatial resolutions in both the atmosphere and troposphere, in order to determine the Earth's atmospheric composition and have a long global ozone record. The instrument covers the range of 270-500 nm with a spectral resolution of about 0.5 nm and is composed of two CCD detectors in the UV (UV1: 270-310 nm and UV2: 310-380 nm) and VIS (350-500 nm) channels. The OMI slit functions were determined from preflight characterization for each spectral pixel \citep{Sun2017a}. The ISRFs of the OMI spectrometer can be found for the UV-2 and VIS bands in the data product of the KNMI Projects made in 2014 in KNMI, De Bilt, The Netherlands. This study focusses on the VIS band without loss of generality. As described in the data product (https://www.knmiprojects.nl/projects/ozone-monitoring-instrument/data-products), the VIS band has only 750 illuminated pixels at the wavelengths ranging from 350 nm to 500 nm. The slit functions of the spectrometer are obtained by averaging the rows 17 to 47 and are derived from the parameterization described in RP-OMIE-KNMI-704. Spline interpolation was used to generate $N_{\lambda}=750$ ISRFs with a sample size $N = 299$. More information on the instrument calibration can be found in \citep{Dobber2006} and in-flight performance was assessed in \citep{Schenkeveld2017}. 

\subsection{TROPOMI}
The Tropospheric Monitoring Instrument (TROPOMI) was launched on the Sentinel-5 Precursor Satellite in October 2017. The instrument is part of the new generation of atmospheric sounding instruments and was designed to continue GOME-2, SCIAMACHY and OMI missions with higher spatial and spectral resolutions \citep{Kleipool2018}. The satellite has four high resolution spectrometers with their own optics and detectors and four different wavelength bands, ranging in the UV (270-320 nm), UVIS (320-500 nm), NIR (675-775 nm) and SWIR (2305-2385 nm) bands with respective resolutions (0.45-0.5 nm), (0.45-0.65 nm), (0.45-0.35 nm) and (0.227-0.225 nm). The ISRFs were first determined during extensive on-ground calibration campaign \citep{vanHees2018} and were then expressed as a convolution between a normal distribution and a uniform distribution. The ISRFs considered in this work were extracted from the TROPOMI Calibration Key Data files whose third version was created in 2018 by KNMI for UV-VIS-NIR spectral bands and in 2016 by SRON for the SWIR band. More details on the data are described in \citep{Tropomi_dataUser}. This study concentrates on the data from the UVIS band having 994 spectral pixels with a spectral range (305-499 nm). Spline interpolation was used to generate $N_{\lambda}=994$ ISRFs with a sample size $N = 257$. More information about the spectrometer can be found in \citep{Tropomi_document}. 

\subsection{OCO-2}
The Orbiting Carbone Observatory 2 (OCO-2) is a NASA Earth observing satellite mission that was launched in July 2014. It is dedicated to the study of atmospheric carbon dioxyde and oxygen and aims at characterizing the global $CO_2$ seasonal cycles and to quantify the sources and sinks of carbone. OCO-2 is composed of three high spectral resolution imaging spectrometers for narrow spectral ranges. The first band is an $O_2$ A-band (757-775 nm) with a spectral resolution of 0.04 nm that is used to constrain the dry-air column abundance and the atmospheric optical path length. The other bands are the weak $CO_2$ band (1594-1627 nm) and the strong $CO_2$ band (2043-2087 nm) with respective spectral resolutions of 0.08 nm and 0.10 nm. These two bands provide information about the $CO_2$ column abundance and aerosol properties. The characterization of ISRFs is highly challenging and crucial due to this high spectral resolution. The ISRFs are measured for each pixel using a tunable diode laser during pre-flight calibration \citep{Lee2017},  and the results are stored in a look-up table. 
The data can be downloaded on the NASA data website EarthDATA. The product considered in this study is the OCO-2 Level 1B Version 11r for science acquired in March 2023 and the fourth footprint is used. Specification on the data product can be found in \citep{OCO2_document}. Some of the ISRFs are declared as unvalid due to radiometric, spatial, spectral or polarization problems (and are thus not considered for ISRF estimation). The ISRFs used for the experiments come from the $O_2$A-band of OCO-2. They were interpolated using splines to generate $N_{\lambda}=859$  ISRFs with a sample size $N = 895$.  More information on OCO-2 spectral calibration and ILS characterization are available in \citet{Lee2017} and \citet{Sun2017b}. 

\subsection{MicroCarb}
MicroCarb is a mission developed by the Centre National des \'Etudes Spatiales (CNES) whose aim is to ensure continuity with other carbone measuring missions such as OCO-2 and GoSat, in order to monitor $CO_2$ fluxes at the Earth surface and determine $CO_2$ atmospheric concentrations. The MicroCarb mission uses a compact and low cost space instrument that will be smaller than the current spectrometers. The instrument is capable of acquiring four spectral bands with a single detector. The first band B1 (758.3-768.3) nm is an $O_2$ band with a spectral resolution of about 0.01 nm. The bands B2 (1596.7-1618.9 nm)  and B3(2023-2051 nm) with respective spectral resolutions of about 0.02 nm and 0.03 nm are sensitive to the concentration of $CO_2$ and have $CO_2$ absorption lines. The last band B4 (1264-1282.2 nm) is a second $O_2$ band with spectral resolution of about 0.02 nm. The wavelengths associated with this last band are closer to the $CO_2$ wavelength and can be used for validation of space-based greenhouse gas observation \citep{Bertaux2020}. The whole dataset has been delivered by the French Space Agency (CNES, Toulouse) containing 1024 ISRFs associated with 1024 spectral measurements for the different bands. The data used for this experiment is the first band of MicroCarb with $N_{\lambda}=1024$ ISRFs and a sample size $N = 895$. More details about MicroCarb can be found on the CNES website (https://microcarb.cnes.fr/en). A particularity of this mission is that the shapes of the ISRFs are strongly dependent on the scene observed by the instrument, which will be discussed in Sect. \ref{Subsubsec:robustnesstoISRFchanges}.

\section{Results and discussion}
\label{section:Results_and_discussion}

\subsection{Numerical experiments and performance evaluation}
The performance of the different ISRF estimation methods is evaluated in terms of ISRF estimation quality and residual between the spectral measurements and their estimates. The quality of ISRF estimation can be quantified by the normalized absolute error between the ISRF and its estimate:
$$
\text{E}_l = \frac{\sum_{n= -N/2} ^{N/2} | I_l(n \Delta) - \hat{I_l} ( n \Delta)|}{\sum_{n= -N/2} ^{N/2} I_l(n \Delta)}.
$$
The residual between the spectral measurements and their estimates is defined by:
$$
\rho = \frac{1}{N_{\lambda}}\sum_{l=1}^{N_{\lambda}} || \boldsymbol{s}_l - \boldsymbol{r}_{l} \hat{\boldsymbol{I}}_l ||_2^2
$$

 In the MicroCarb mission, the ISRFs are considered to be well estimated when their normalized errors satisfy $\text{E}_l <1\%$ for each wavelength. The proposed SPIRIT method is compared to the parametric methods based on Gaussian and Super-Gaussian models. The parameters of these models are estimated using the non-linear least squares algorithm based on the Nelder-Mead optimization algorithm \citep{Nelder_Mead} (MATLAB function \textit{fminsearch}). This iterative algorithm requires an initialisation and a stopping criterion. For the initialization of the Gaussian model, the mean $\mu_{\text{G}_0}$ was set to the sample mean of the ISRFs, the Full Width at Half Maximum (FWHM) was used for the standard deviation $\sigma_{\text{G}_0}$ and the amplitude was initialized as $A_{\text{G}_0} = (2 \pi \sigma_{\text{G}_0})^{-1/2}$. For the Super-Gaussian model, the initialization was defined as $\mu_{\text{SG}_0} = \mu_{\text{G}_0}$,  $k_{\text{SG}_0} = 2$, $w_{\text{SG}_0} = \sqrt{2} \sigma_{\text{G}_0}$ and $A_{\text{SG}_0} = \frac{k_{\text{SG}_0}}{2 w_{\text{SG}_0} } \Gamma(1/k_{\text{SG}_0})$ where $\Gamma$ is the gamma function. The algorithm was stopped after a maximum number of iterations equal to $20000$. The proposed SPIRIT method requires to build a dictionary of representative spectral responses. This dictionary was constructed using an SVD of ISRF examples or estimated using the K-SVD algorithm.  Two different sparse coding methods based on LASSO \citep{tibshirani96regression}  and OMP are investigated after dictionary construction. The first method uses a MATLAB implementation of LASSO with a parameter $\mu > 0$ adjusted to obtain a desired  number of atoms. The non zero coefficients obtained with LASSO were re-estimated in order to reduce the shrinking bias of LASSO \citep{LASSO_LS}. All the ISRF estimation methods were evaluated for the six instruments presented in Sect. \ref{section:Instruments_and_datasets}.   

 \subsection{ISRF estimation for a given wavelength}
\begin{figure}[htb]
  
 \centering
      a)
            \includegraphics[trim={0cm 0cm 0cm 0cm}, clip,width=7cm]{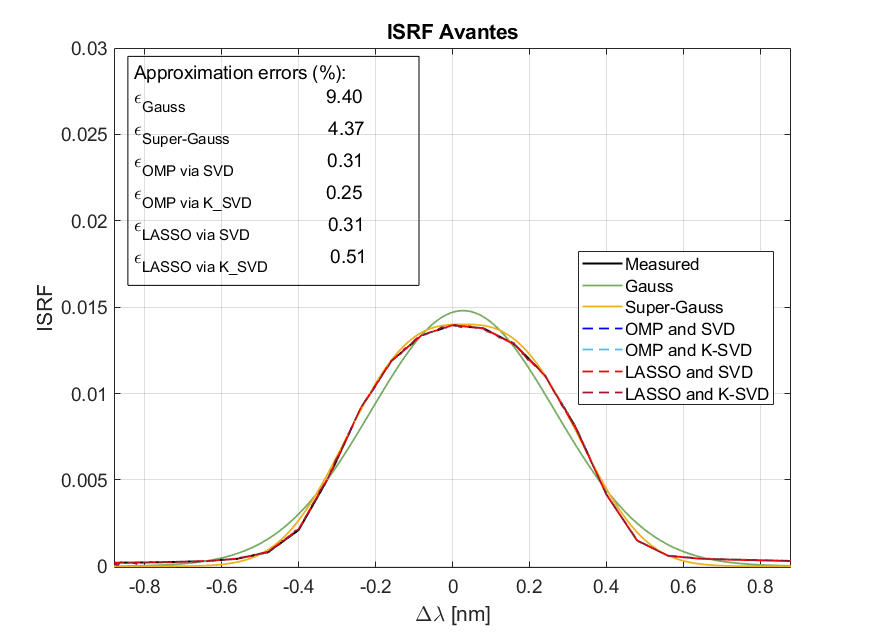}
        b)
            \includegraphics[trim={0cm 0cm 0cm 0cm}, clip,width=7cm]{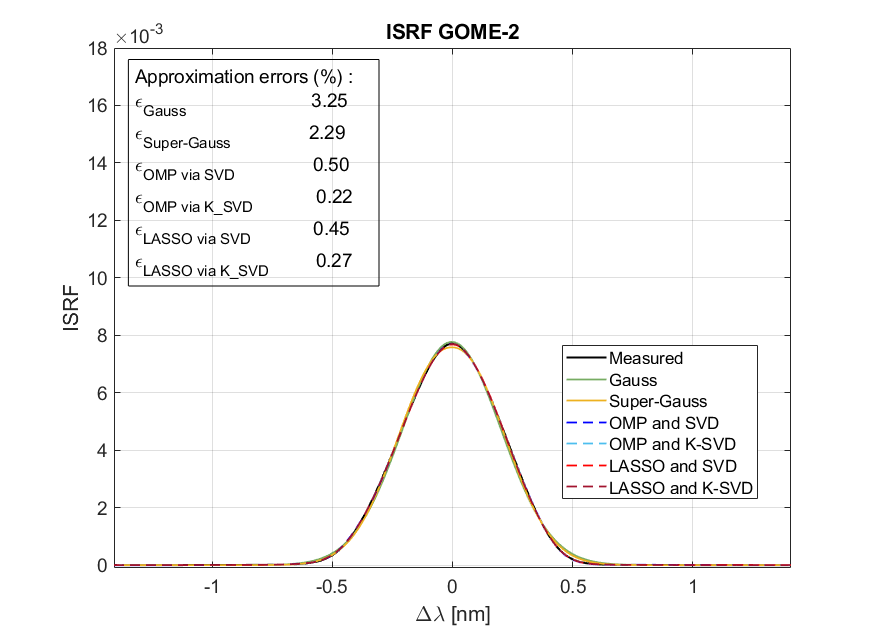}
\\
       c)    
            \includegraphics[trim={0cm 0cm 0cm 0cm}, clip,width=7cm]{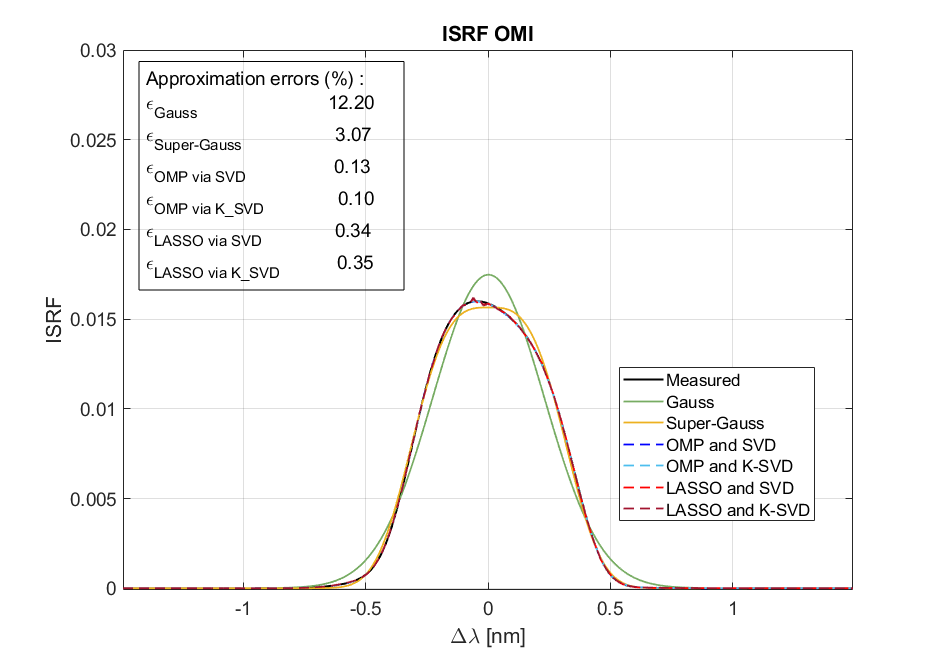}
        d)
            \includegraphics[trim={0cm 0cm 0cm 0cm}, clip,width=7cm]{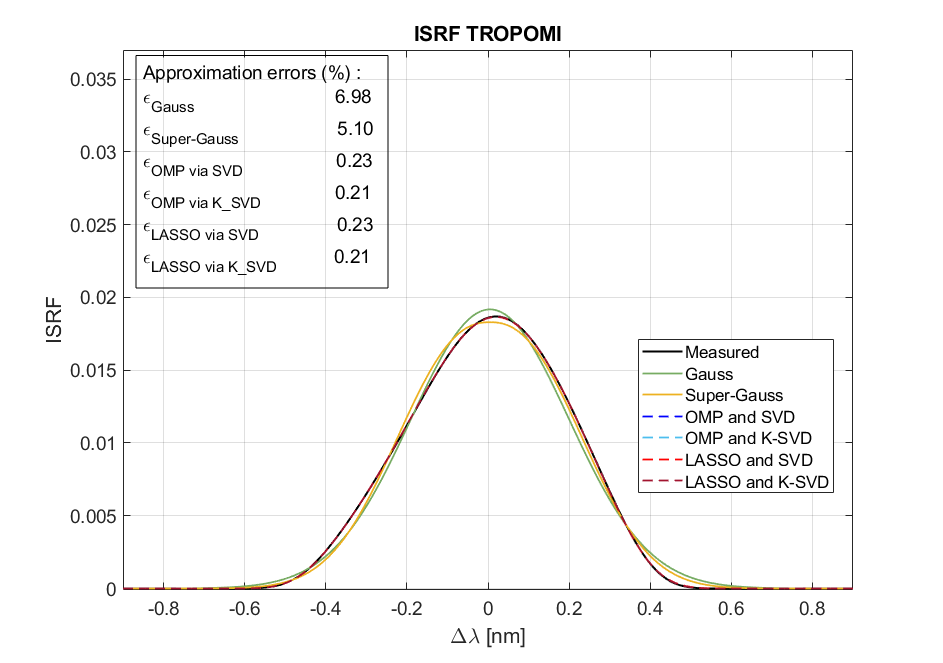}
\\
       e)
            \includegraphics[trim={0cm 0cm 0cm 0cm}, clip,width=7cm]{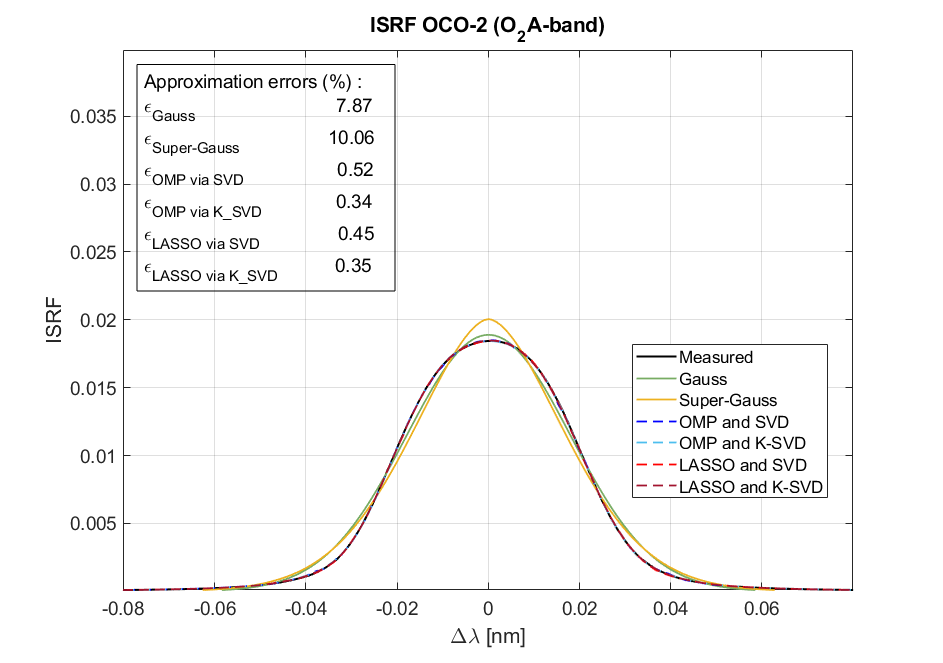}
        f)
            \includegraphics[trim={0cm 0cm 0cm 0cm}, clip,width=7cm]{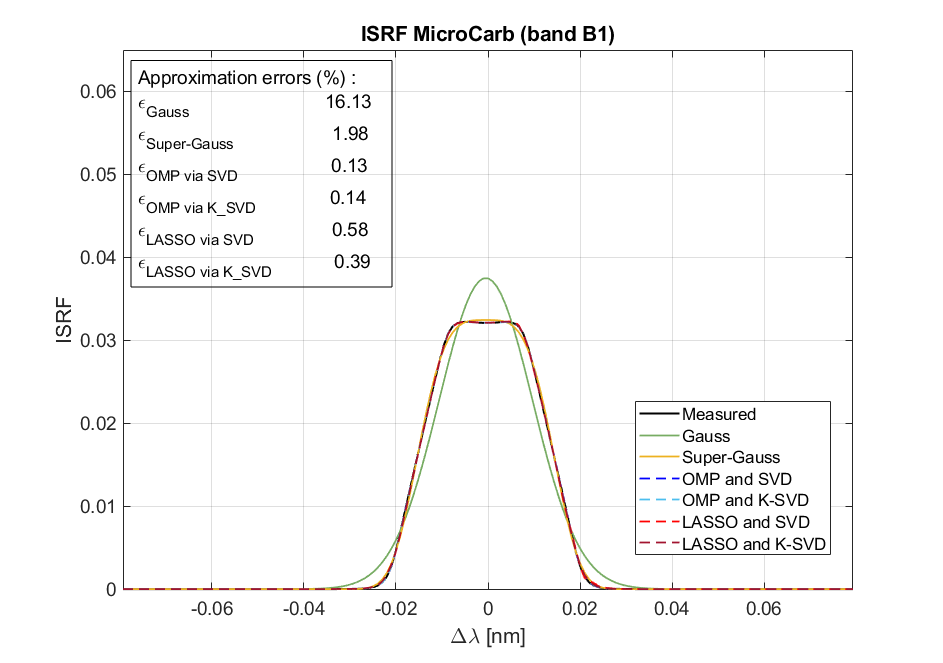}
     \caption{Examples of ISRFs and their estimates using parametric methods and SPIRIT for ISRFs centered at 404.25 nm for Avantes (a), centered at 430 nm for GOME-2 (b), OMI (c) and TROPOMI (d) and centered at 763.01 nm for OCO-2 (e) and MicroCarb (f).}
    \label{figure:illustration_single_ISRF}
\end{figure}

Examples of ISRFs and their estimates are displayed in Fig. \ref{figure:illustration_single_ISRF}. First, it is interesting to note that the ISRF shapes can be significantly different depending on the considered wavelength and the instrument. This observation suggests that the dictionary must be adapted to the spectrometer. The results shown in Fig. \ref{figure:illustration_single_ISRF} clearly illustrate the advantage of using SPIRIT for ISRF estimation, which leads to normalized estimation errors of less than $1\%$, significantly below those obtained using the parametric estimation methods. 
A comparison between the different sparse approximations (OMP, LASSO) and dictionaries (SVD, K-SVD) that can be used by SPIRIT shows that OMP works better than LASSO for this example. Also, with the exception of the GOME-2 and the $O_2$A-band of OCO-2 spectrometers, using the K-SVD algorithm does not significantly improve the results with respect to SVD, at the price of a higher computational complexity. A preliminary conclusion is that the use of SPIRIT with a combination of OMP and SVD will be the method of choice.

\begin{figure}[htb] 
\begin{minipage}[t]{0.5\linewidth}\centering
            \includegraphics[trim=50 0 50 0, clip,width=\linewidth]{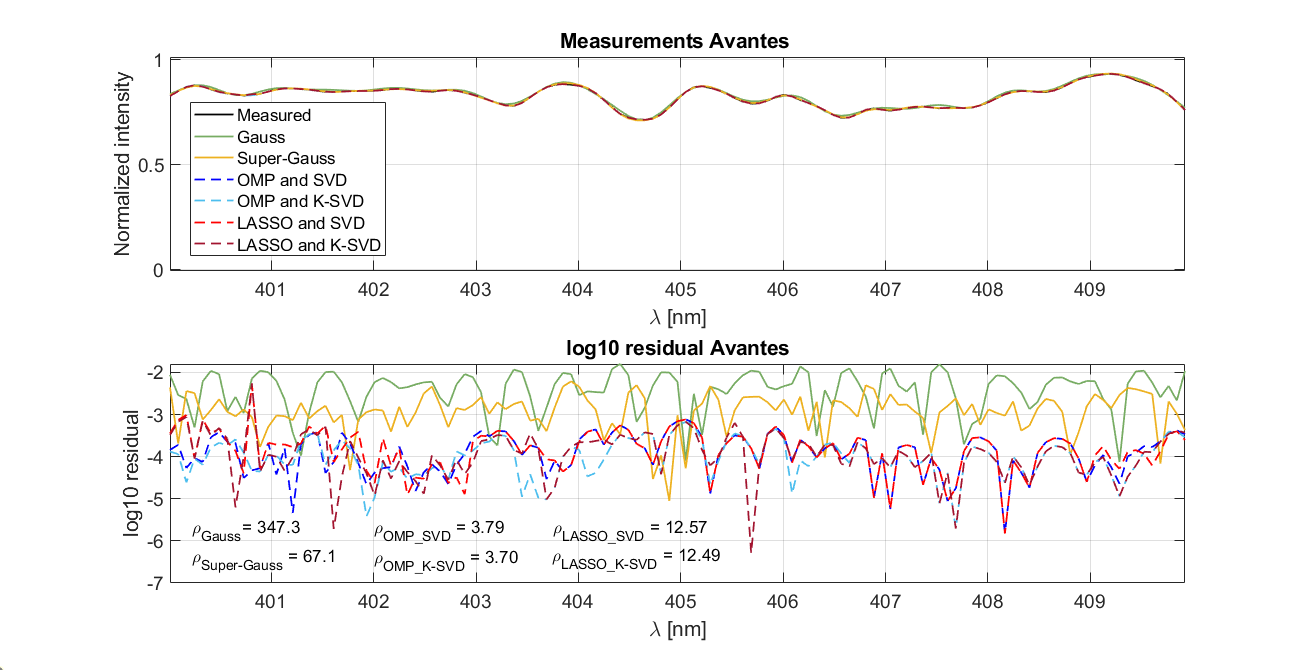}
            \includegraphics[trim=50 0 50 0, clip,width=\linewidth]{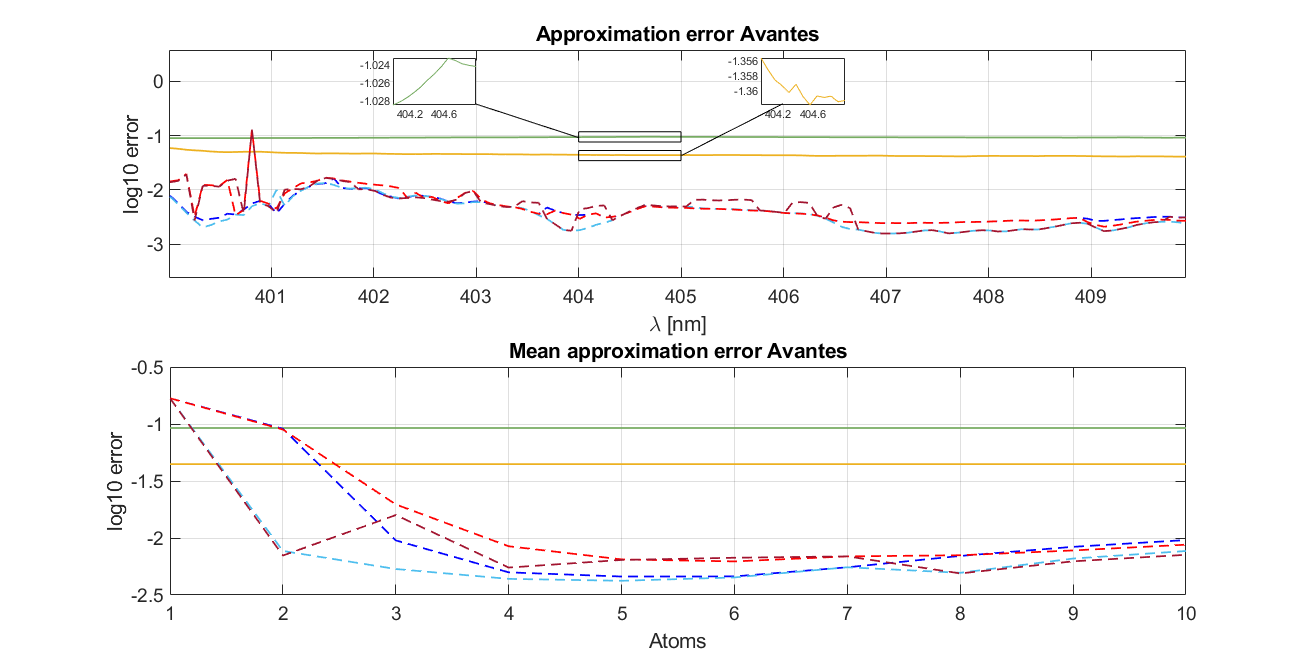}\\
            (a)
\end{minipage}%
\begin{minipage}[t]{0.5\linewidth}\centering
            \includegraphics[trim=50 0 50 0, clip,width=\linewidth]{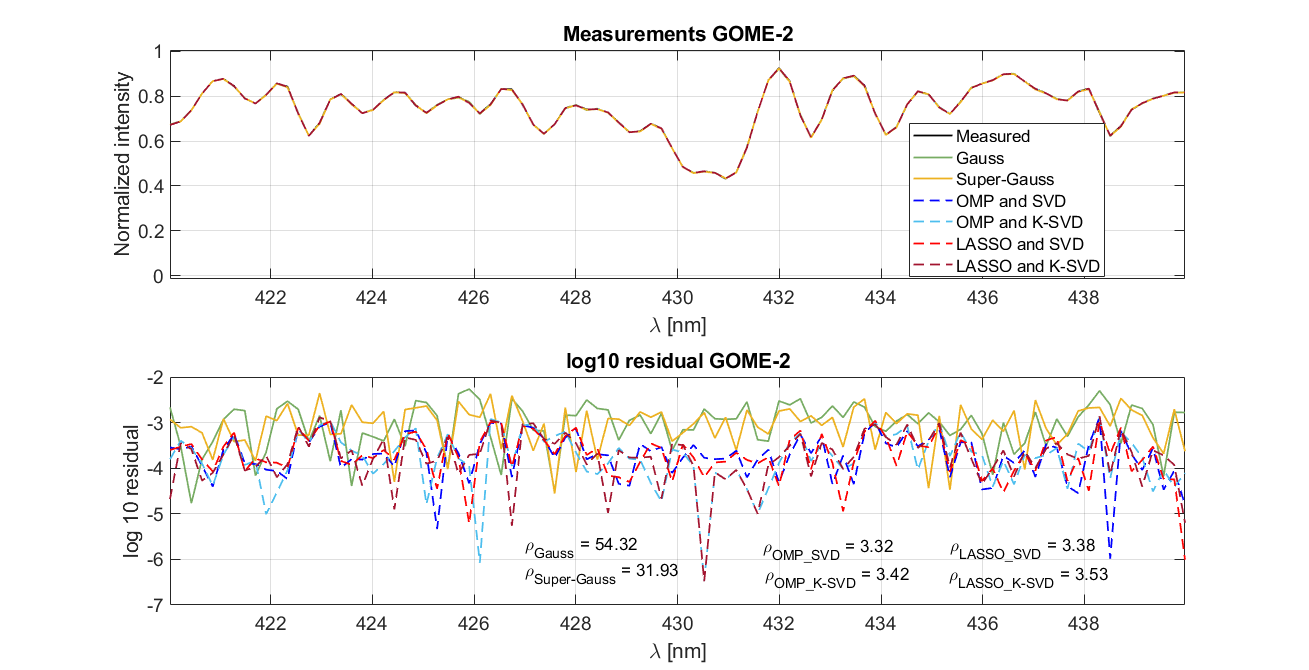}
            \includegraphics[trim=50 0 50 0, clip,width=\linewidth]{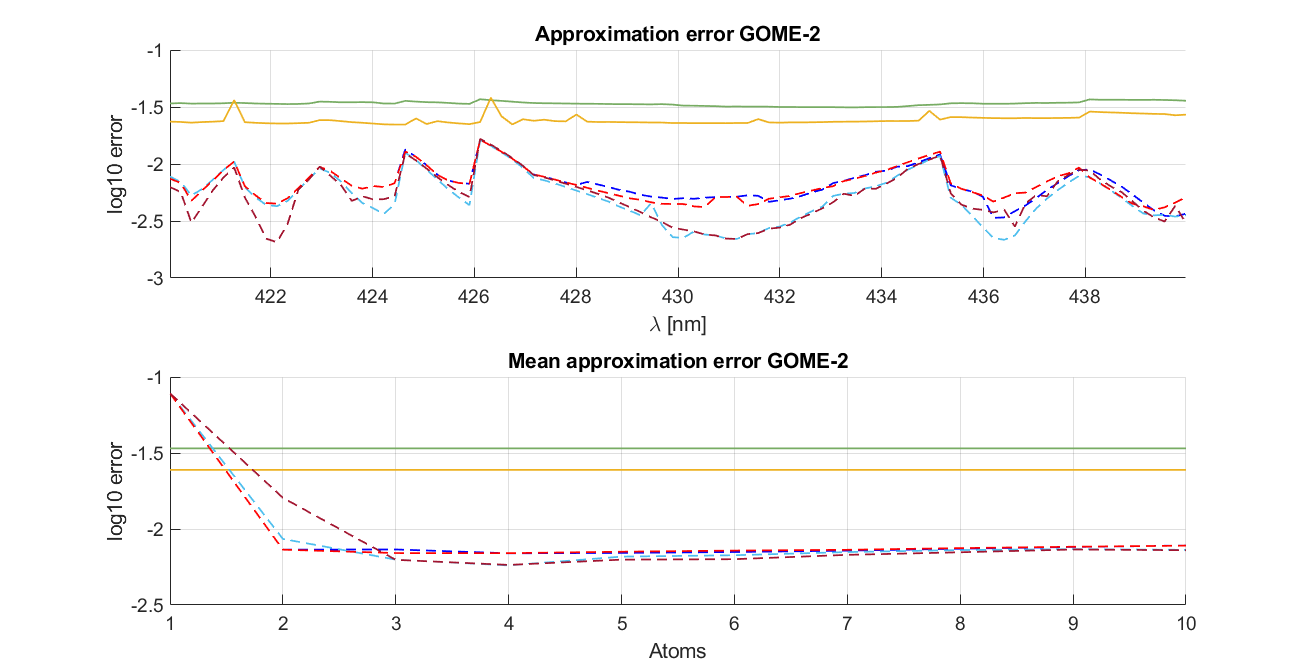}\\
            (b)
\end{minipage}
    \caption{Avantes (a) and GOME-2 (b) ISRF estimates using different methods (Gauss, Super-Gauss, OMP and LASSO, SVD and K-SVD).}
    \label{figure:Avantes_GOME-2_spectrometer}
\end{figure}

\begin{figure}[htb] 
\begin{minipage}[t]{0.5\linewidth}\centering
            \includegraphics[trim=50 0 50 0, clip,width=\linewidth]{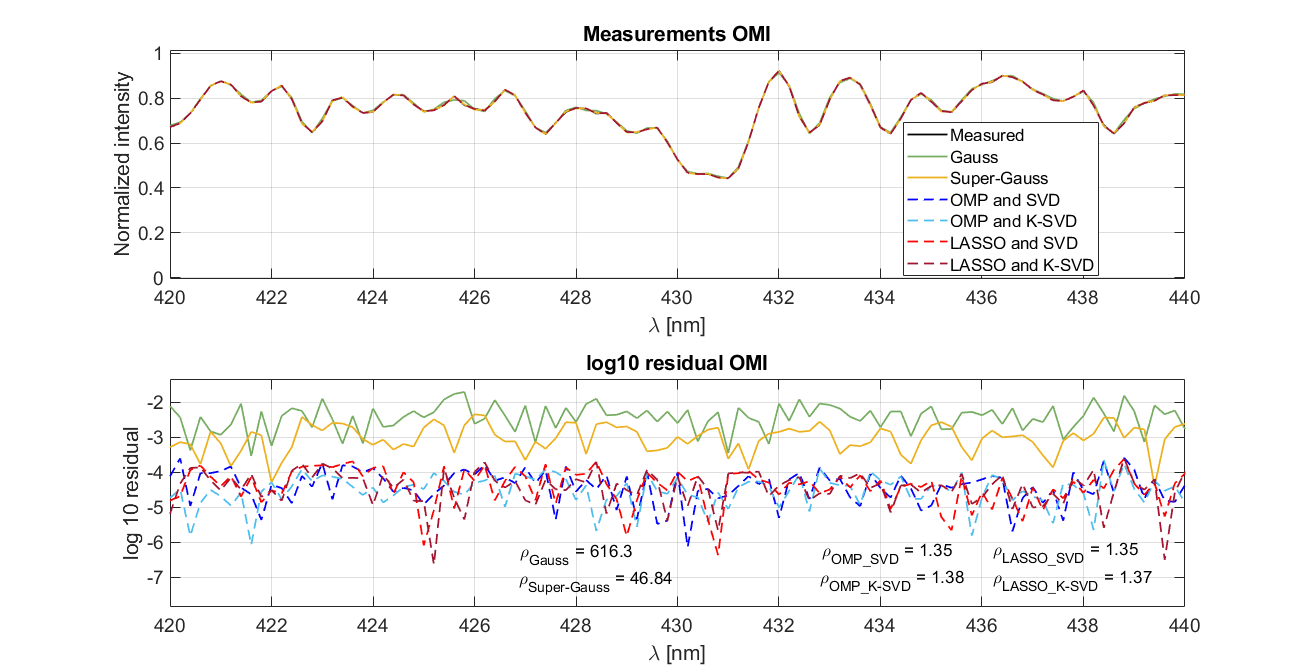}
            \includegraphics[trim=50 0 50 0, clip,width=\linewidth]{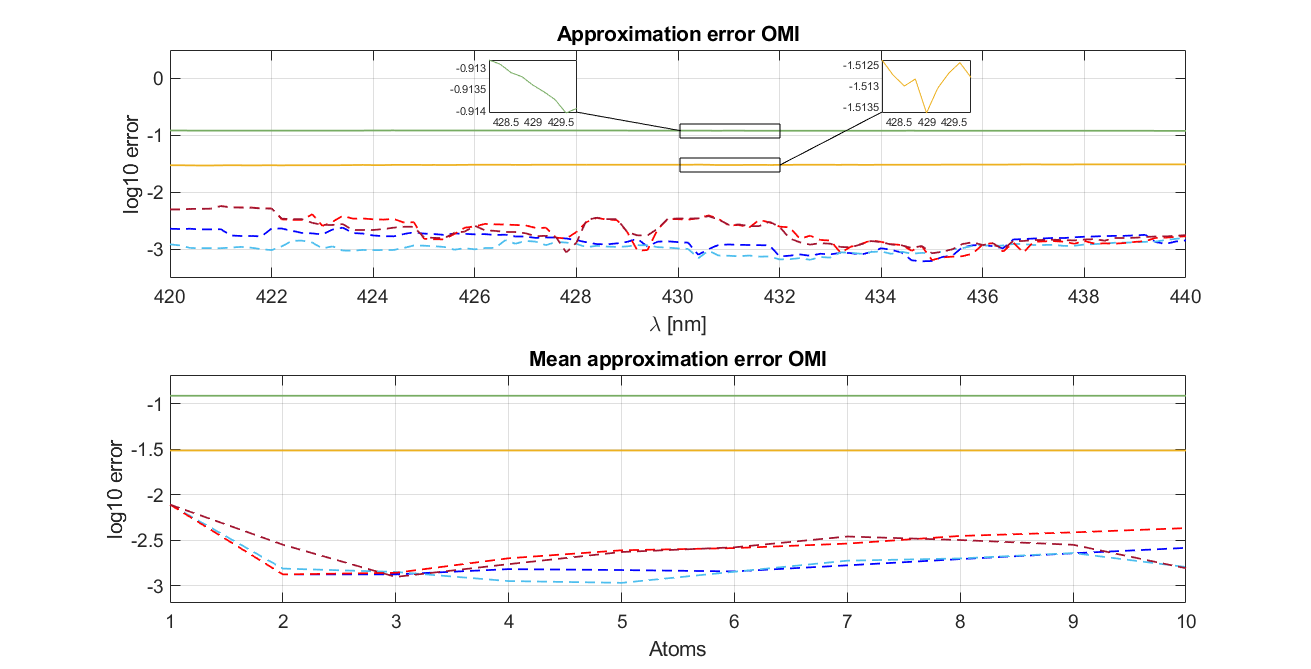}\\
            (a)
\end{minipage}%
\begin{minipage}[t]{0.5\linewidth}\centering
            \includegraphics[trim=50 0 50 0, clip,width=\linewidth]{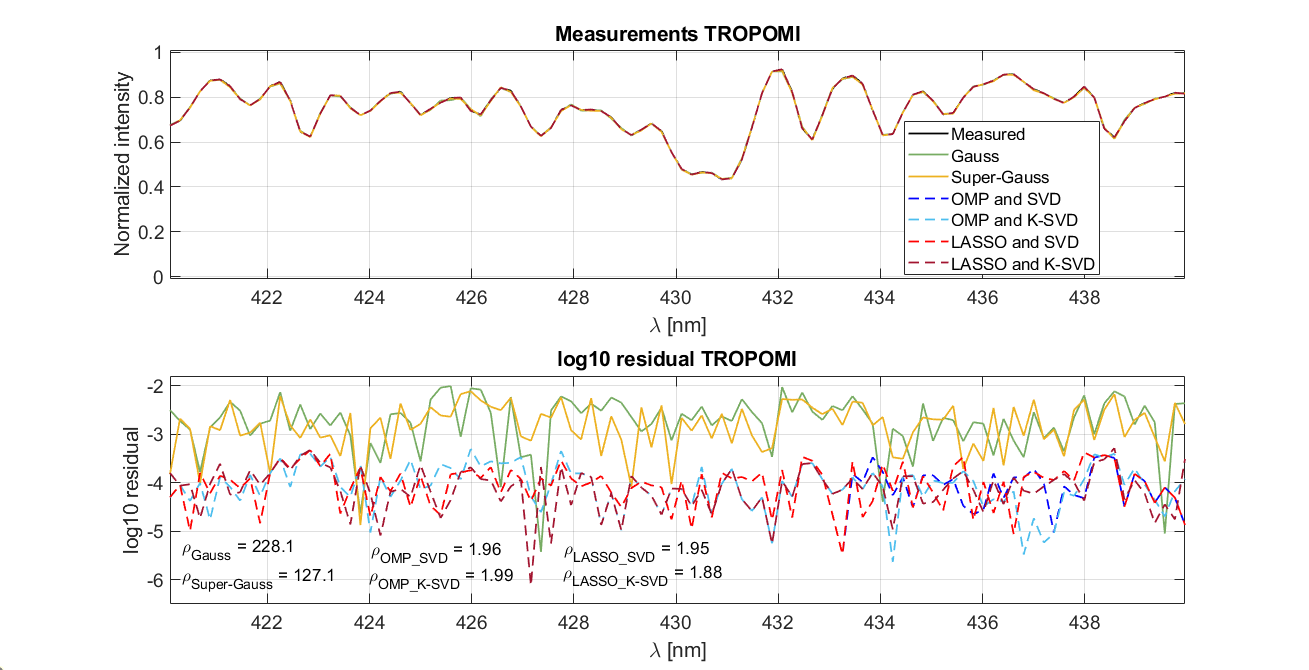}
            \includegraphics[trim=50 0 50 0, clip,width=\linewidth]{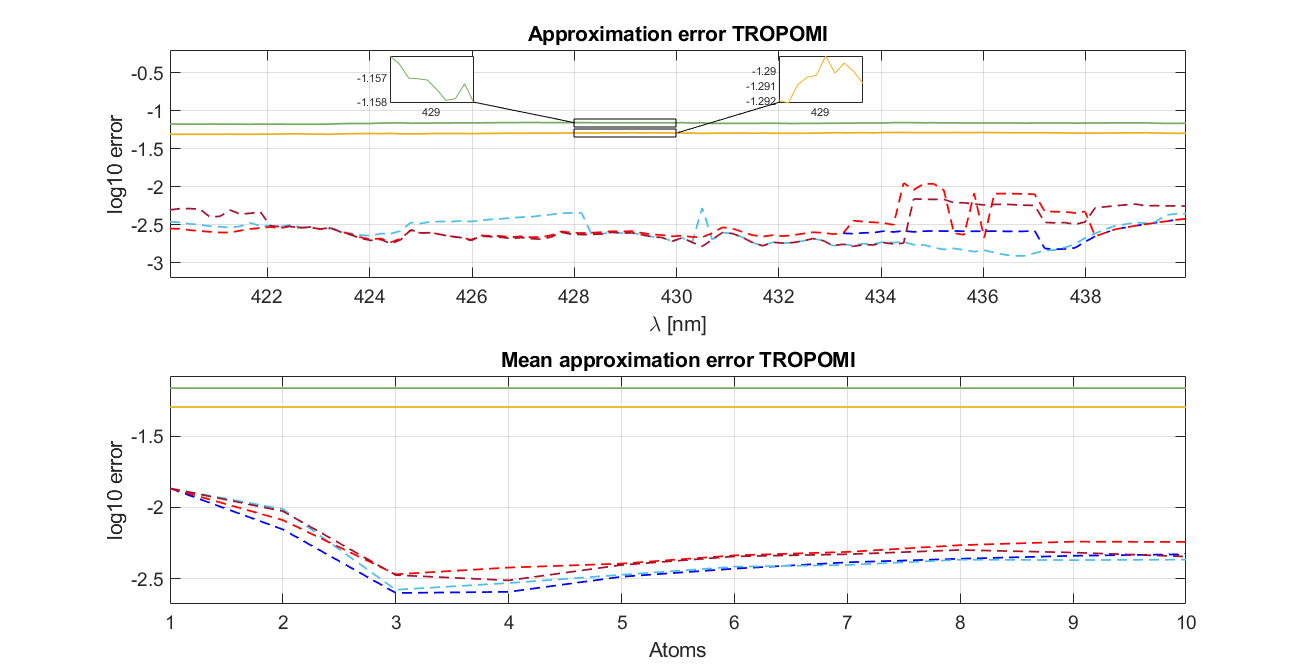}\\
            (b)
\end{minipage}
    \caption{OMI (a) and TROPOMI (b) ISRF estimates using different methods (Gauss, Super-Gauss, OMP and LASSO, SVD and K-SVD).}
    \label{figure:OMI_TROPOMI_spectrometer}
\end{figure}

\begin{figure}[htb] 
\begin{minipage}[t]{0.5\linewidth}\centering
            \includegraphics[trim=50 0 50 0, clip,width=\linewidth]{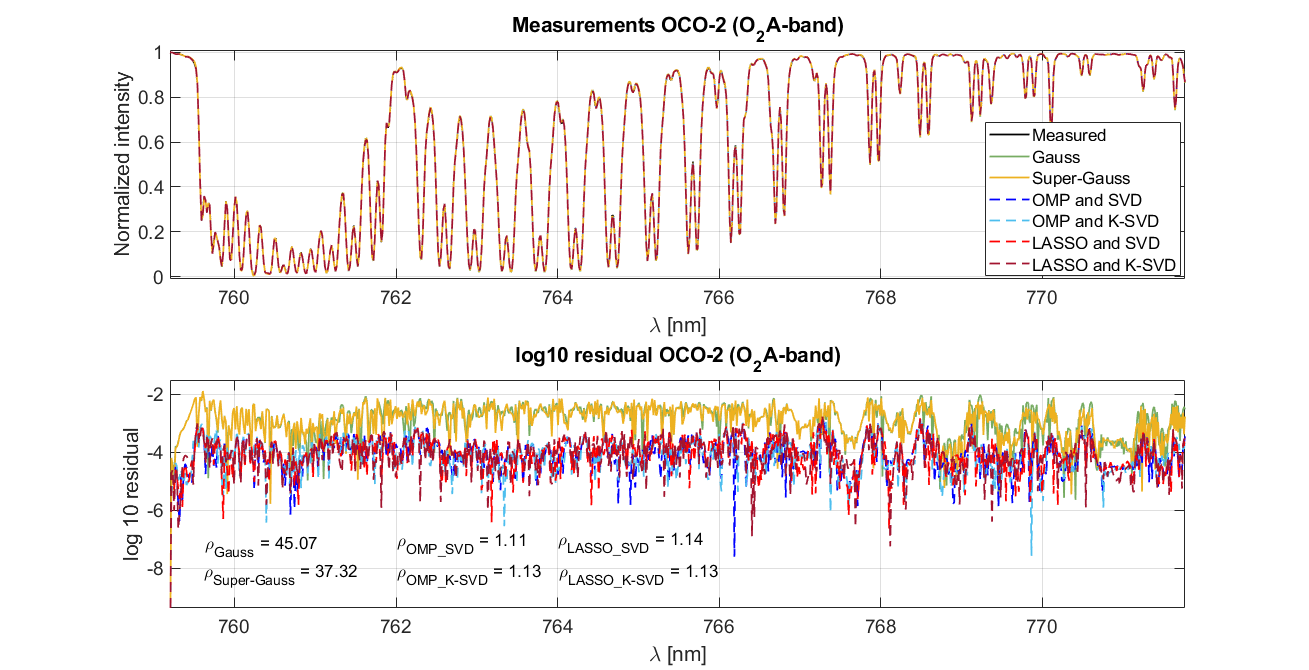}
            \includegraphics[trim=50 0 50 0, clip,width=\linewidth]{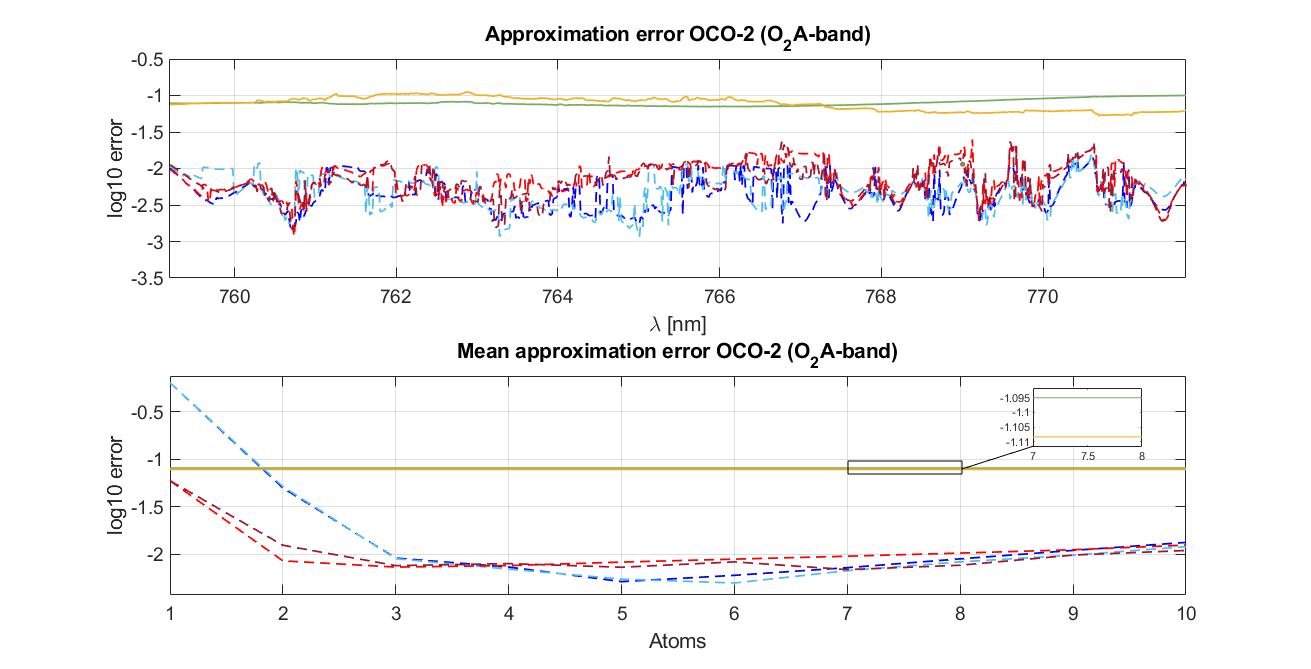}\\
            (a)
\end{minipage}%
\begin{minipage}[t]{0.5\linewidth}\centering
            \includegraphics[trim=50 0 50 0, clip,width=\linewidth]{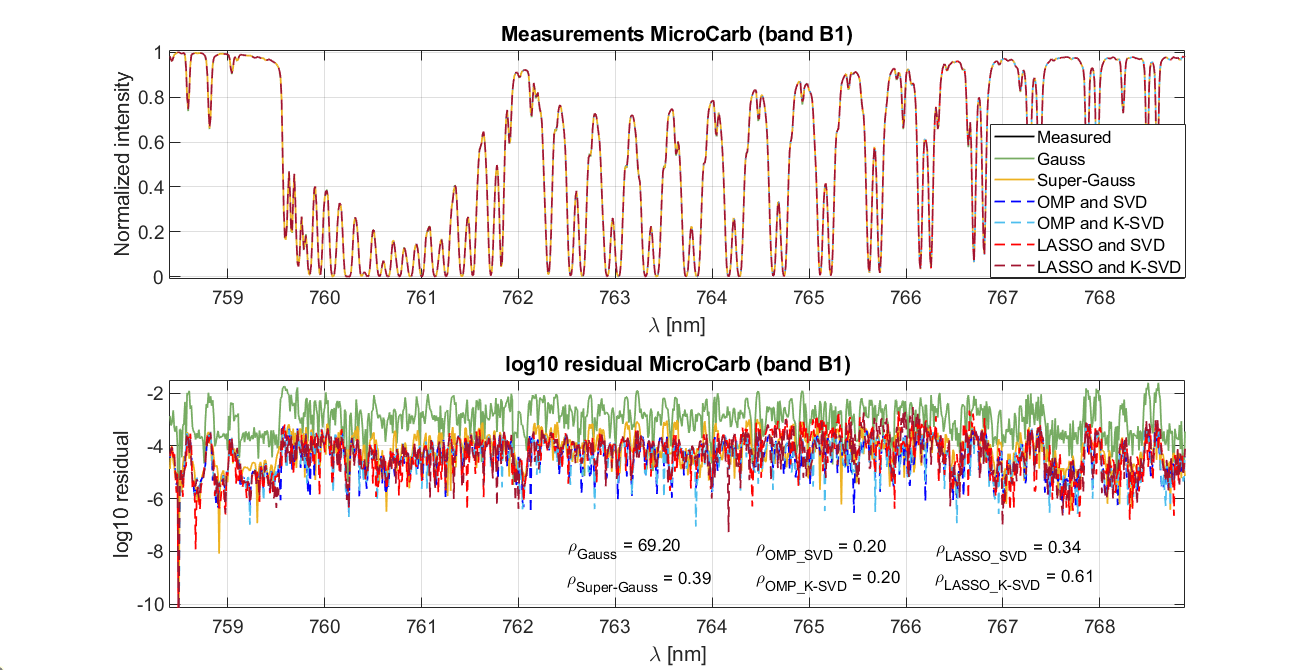}
            \includegraphics[trim=50 0 50 0, clip,width=\linewidth]{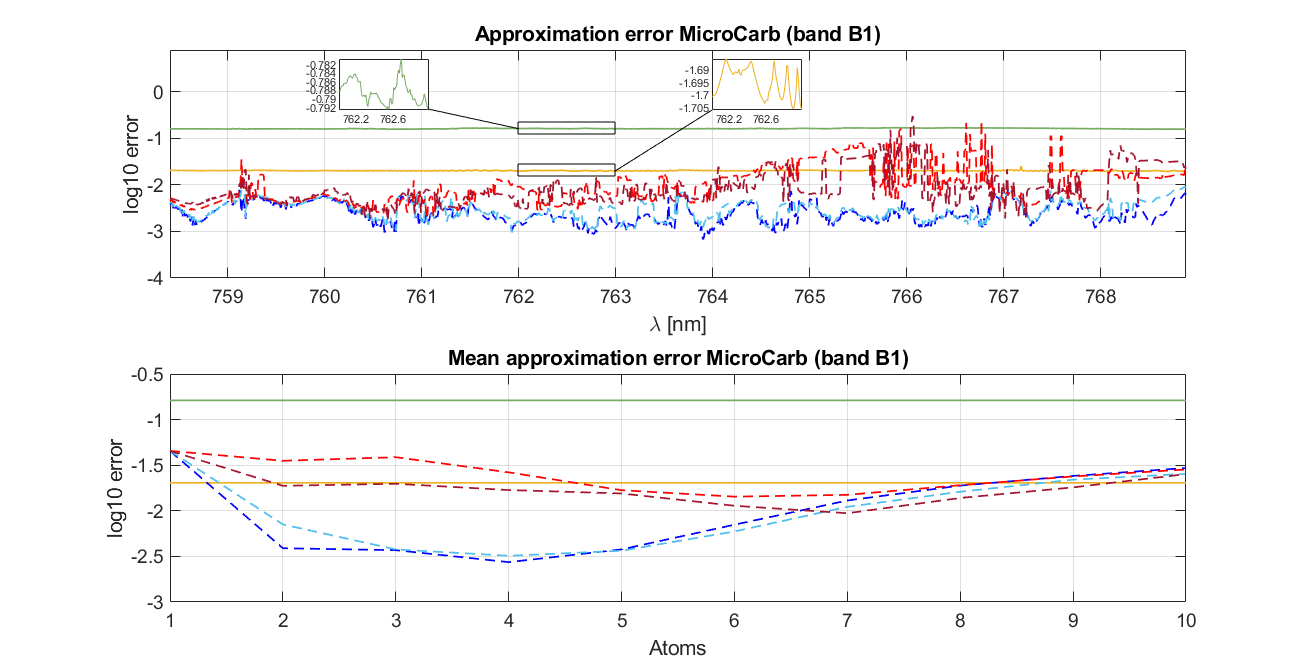}\\
            (b)
\end{minipage}
    \caption{OCO-2, $O_2$ A-band (a) and MicroCarb, band 1 (b) ISRF estimates using different methods (Gauss, Super-Gauss, OMP and LASSO with SVD or K-SVD).}
    \label{figure:OCO_2_CNES_spectrometer}
\end{figure}

 \subsection{Performance versus wavelength}
 \label{sec:perfwav}
Figs.~\ref{figure:Avantes_GOME-2_spectrometer} to \ref{figure:OCO_2_CNES_spectrometer} display performance results for all wavelengths: spectral measurements and their reconstructions using the different methods (top rows), absolute differences between the normalized spectral measurements and their normalized approximations (in logarithmic scale, second rows) and the associated residuals $\rho$, comparison between the ISRF approximation errors for different wavelengths (third rows), ISRF errors averaged over all the wavelengths as a function of the number of atoms selected from the dictionary (bottom rows). The number of atoms for the first three rows is determined as the value of $K$ minimizing the average ISRF error reported in the bottom row.
  
\subsubsection{Avantes}

The measurements considered in this section correspond to the wavelengths ranging from 400 to 410 nm, to make a fair comparison with the results of \citep{Beirle2017}. The bottom plot of Fig.~\ref{figure:Avantes_GOME-2_spectrometer} (a)  shows that $K = 5$ atoms is optimal for this instrument. Approximation errors for ISRF estimates using  the Super-Gauss parameterization are lower than those obtained using the Gaussian model, which confirms the results of \citep{Beirle2017}. The proposed SPIRIT method allows us to obtain normalized errors less than $1\%$, which are significantly below those of parametric methods. The OMP algorithm tends to provide slightly better results than LASSO for this instrument. Note that the differences obtained using SVD and KSVD for building the dictionary are found to be insignificant.

\subsubsection{GOME-2}
The measurements considered in this paper have wavelengths ranging  from 420 nm to 440 nm. The bottom plot of Fig.~\ref{figure:Avantes_GOME-2_spectrometer} (b) shows that choosing $K = 4$ atoms chosen from the dictionary leads to the best ISRF estimations. Inspection of the other rows of Fig.~\ref{figure:Avantes_GOME-2_spectrometer} (b) indicates again that the Super-Gauss parameterization performs better than the Gaussian one, but significantly worse than the proposed sparse representation. The K-SVD algorithm provides slightly better estimation results than SVD for this dataset, notably in the range 428 nm to 435 nm.

\subsubsection{OMI}
The measurement considered in this section are for wavelengths 420 to 440 nm as for the GOME-2 spectrometer. The measured spectrum is reconstructed with the sparse representation methods for $K = 5$ atoms chosen from the dictionary, which is constructed using SVD or K-SVD.
For this spectrometer, Fig.~\ref{figure:OMI_TROPOMI_spectrometer} (a) highlights particularly well the importance of choosing a method other than the Gaussian model. The Super-Gaussian parameterization delivers better results in terms of measurement fit and ISRF approximation. However, methods based on sparse representation provide much better results, achieving mean approximation errors as low as $0.1\%$.

\subsubsection{TROPOMI}
Fig.~\ref{figure:OMI_TROPOMI_spectrometer} (b) shows that sparse methods applied to TROPOMI data require a smaller number of atoms to minimize both the measurement fitting errors and the ISRF approximation errors. The proposed SPIRIT method yields significantly better performance than state-of-the-art parametric models, leading to more than one magnitude smaller ISRF approximation errors.

\subsubsection{OCO-2}
The OCO-2 measurements considered in this section are obtained using the data for the $O_2$ band (757-772 nm). The measured spectrum is reconstructed with the sparse representation methods for K = 5 atoms chosen from the dictionary (constructed using SVD or K-SVD). Results shown in Fig.~\ref{figure:OCO_2_CNES_spectrometer} (a) indicate that the Super-Gauss model delivers slightly better results than the Gaussian model in terms of residual error and mean ISRF approximation error. However, for the low wavelengths of the band, the ISRF approximation errors are slightly higher with the Super-Gaussian model, possibly due to the shapes of the ISRFs that are closer to be Gaussian, or to the spectral sampling position which can induce some biases for this spectrometer \citep{Sun2017b}. Both parametric models yield close to $10\%$ ISRF approximation errors.
The proposed sparse representation approach again yields far better ISRF approximations and measurement fits, with the best results obtained using OMP and SVD.

\subsubsection{MicroCarb}
The measurements used in this part were simulated by the CNES for the B1 wavelength range (758.4-768.9 nm). Results displayed in Fig.~\ref{figure:OCO_2_CNES_spectrometer} (b) show that, for the MicroCarb spectrometer, the use of the Super-Gauss parameterization  reduces the measurement fitting and ISRF approximation errors compared with the Gaussian model. However, the methods based on sparse representations are performing significantly better and reach results below the $1\%$ error criterion for this mission. In the case of this spectrometer, LASSO has more difficulty to approximate the ISRFs than OMP, possibly due to the $l_1$ regularization, LASSO may struggle to select the correct coefficients. Thus, sparse coding should be carried out using OMP, which minimizes the approximation errors to a much greater extent, enabling errors of the order of $0.1\%$ to be achieved for certain wavelengths.

\subsubsection{Discussion}

Overall, the conclusions from these experiments are as follows. First, the Super-Gaussian parameterization generally yields better performance than the Gaussian one, corroborating the results reported in \citep{Beirle2017}. However, the normalized ISRF approximation errors obtained with these parametric methods are consistently larger than $1\%$, for most instruments and wavelengths, contrary to the proposed SPIRIT approach based on sparse approximations of ISRFs in a suitable dictionary. 
This result is due to the fact that the ISRF shapes depend strongly on the spectrometer and can vary across wavelengths, which cannot be accommodated easily with a simple parametric model. On the contrary, decompositions in appropriate dictionaries that depend on the spectrometer and the chosen wavelength band offer sufficient flexibility for all use cases considered in this paper. 
Regarding the estimation algorithms, SVD provides an estimation performance close to K-SVD except in some very specific cases. LASSO is usually outperformed by OMP. These results overall suggest the use of SVD for building the dictionary and OMP for ISRF estimation.

\subsection{Robustness analysis \& ablation study}
\subsubsection{Robustness to noise}
To study the robustness of the different ISRF estimation methods to the presence noise, white Gaussian noise was added to the spectral measurements with several signal to noise ratio (SNR) levels above and below the value of $55$dB that is expected for the MicroCarb mission. Table \ref{Table_SNR} reports the obtained residual approximation errors and the normalized average ISRF approximation errors for all instruments. 
Approximation errors less than $< 1\%$ are highlighted in blue. 
These results show that the proposed sparse representations yield normalized errors below $1\%$ for SNRs larger than $20$ dB. Note that the use of LASSO for sparse coding leads to larger errors for the MicroCarb spectrometer, again leading to the recommandation of using OMP instead of LASSO.
The proposed ISRF estimation methods based on sparse representations globally yield the smallest approximation and residual errors, hence the best estimation results. Note that the errors obtained with the parametric models do not vary significantly with the noise level, except for the smallest SNR, indicating that errors due to model misfit are larger than those induced by the noise degradations. To conclude, OMP combined with SVD provides the overall best results for ISRF estimation, also in the presence of additive noise.

\begin{table}[htbp!]
      \caption{Mean residual and approximation errors for different SNRs and different methods (Gauss (G), Super-Gauss (SG), OMP and LASSO, SVD and K-SVD).} 
    \label{Table_SNR}
  
    \begin{center}
    \setlength{\tabcolsep}{1.8pt}
     \small

    \begin{tabular}{|cc||c|c|c|c|c|c||c|c|c|c|c|c|}
               \hline
\multicolumn{2}{|c||}{}&\multicolumn{6}{c||}{Normalized approximation error (\%)}&\multicolumn{6}{c|}{Residual error}\\
        \hline
      Instrument & / SNR   & G & SG & OMP & OMP & LASSO& LASSO& G & SG & OMP& OMP & LASSO & LASSO   \\ 
             &   &  &  &SVD& K-SVD & SVD & K-SVD &  &  & SVD & K-SVD & SVD & K-SVD   \\\hline 
       & 20 dB &  9.27 &  4.78 & \textbf{3.06} & 3.30 & 6.04 &   5.89 &   645.5 &   351.7 &  271.1 &   \textbf{270.3} &   311.1 & 303.2 \\
       & 40 dB &  9.22 &   4.45 &   \color{blue}{\textbf{0.50}} & \color{blue}{0.51} & \color{blue}{0.80} &   \color{blue}{0.83} &  351.2 &  70.06 &  6.50 &  \textbf{6.40}  & 15.21 & 15.13\\
      Avantes & 55 dB &  9.22 &   4.45 &  \color{blue}{0.46} & \color{blue}{\textbf{0.43}} & \color{blue}{0.65} &   \color{blue}{0.65} &   347.6 & 67.23 &  3.87 &   \textbf{3.78} &   12.64 &  12.58 \\   
       & 80 dB &  9.22 &  4.45 &   \color{blue}{0.46} & \color{blue}{\textbf{0.42}} & \color{blue}{0.65} &   \color{blue}{0.64} &  347.3 & 67.12 & 3.79 &   \textbf{3.70} &   12.57 &   12.49 \\
       & 120 dB &   9.22 &   4.45 & \color{blue}{0.46} &   \color{blue}{\textbf{0.42}} & \color{blue}{0.65} &   \color{blue}{0.64} &  347.3 &  67.12 &  3.79 &  \textbf{3.70} & 12.57 &  12.49 \\ \hline
       & 20 dB &  3.47 &   2.76 &   2.63 & 2.54 & 2.61 &  \textbf{2.38} &   342.1 &   321.2 &  \textbf{264.4} &   265.8 &   266.3 & 269.1 \\
       & 40 dB &   3.39 &   2.49 &  \color{blue}{0.72} & \color{blue}{\textbf{0.65}} & \color{blue}{0.74} &   \color{blue}{\textbf{0.65}} &  57.26 &  37.45 &  \textbf{6.04} &   6.08 &   6.08 & 6.25 \\
      GOME-2 & 55 dB &  3.39 &   2.46 &   \color{blue}{0.69} & \color{blue}{0.59} & \color{blue}{0.69} &   \color{blue}{\textbf{0.58}} &    54.42 &  32.79 &   \textbf{3.42} &  3.50 &   3.46 & 3.62  \\   
       & 80 dB & 3.39 &  2.44 &   \color{blue}{0.69} & \color{blue}{\textbf{0.58}} & \color{blue}{0.69} & \color{blue}{\textbf{0.58}} & 54.32 &  30.33 &   \textbf{3.33} &   3.42 &   3.38 & 3.53  \\
      & 120 dB &  3.39 &   2.45 &   \color{blue}{0.69} & \color{blue}{\textbf{0.58}} & \color{blue}{0.69} &   \color{blue}{\textbf{0.58}} & 54.32 &  32.00 &   \textbf{3.32} &  3.42 &  3.38 &   3.53  \\ \hline
       & 20 dB & 12.18 &   3.30 &   \textbf{2.71} & 2.95 & 6.48 &   6.68 & 904.6 &   329.3 &   261.4 &   265.6 &  \textbf{255.5} & 259.3  \\
       & 40 dB & 12.17 &   3.06 &   \color{blue}{0.30} & \color{blue}{\textbf{0.28}} & \color{blue}{0.68} & \color{blue}{0.71} &   619.0 &   49.39 &  3.97   & 4.04 &  \textbf{3.94} &   3.97 \\
      OMI & 55 dB & 12.17 &  3.06 &   \color{blue}{0.15} & \color{blue}{\textbf{0.12}} & \color{blue}{0.26} &  \color{blue}{0.26} &  616.4 &  46.88 &   \textbf{1.43} &   1.46 &  \textbf{1.43} &   1.45 \\   
       & 80 dB &12.17 &   3.06 &  \color{blue}{0.15} & \color{blue}{\textbf{0.11}} & \color{blue}{0.25} & \color{blue}{0.24} &  616.3 & 46.84 &  \textbf{1.35} &   1.38 &   \textbf{1.35} &   1.37
 \\ 
       & 120 dB &  12.17 &  3.06 &   \color{blue}{0.15} & \color{blue}{\textbf{0.11}} & \color{blue}{0.25} & \color{blue}{0.24} &  616.3 &   46.84 &  \textbf{1.35} & 1.38 &   \textbf{1.35} & 1.37   \\ \hline
       & 20 dB &6.88 &   5.16 &  \textbf{2.16} &   2.17 &   3.39 & 3.42 & 517.6 &  410.9 &  \textbf{273.9} &   275.1 &   274.2   & 276.5  \\
       & 40 dB & 6.87 &   5.07 &   \color{blue}{\textbf{0.30}} &   \color{blue}{0.34} &   \color{blue}{0.41} & \color{blue}{0.38} &  230.9 &   129.7 &  4.70 &   4.74 &  4.69 &  \textbf{4.61} \\
      TROPOMI & 55 dB & 6.87 &   5.07 &  \color{blue}{ \textbf{0.25} } &   \color{blue}{0.26} &   \color{blue}{0.33} & \color{blue}{0.33} & 228.2 &  127.1 &   2.04 &  2.05 & 2.02 &  \textbf{1.96}  \\   
       & 80 dB &  6.87 &   5.07 &   \color{blue}{\textbf{0.25}} &   \color{blue}{0.26} & \color{blue}{0.34} & \color{blue}{0.33}&   228.1 & 127.1  & 1.96 &   1.99 &   1.95 &   \textbf{1.88} \\ 
       & 120 dB & 6.87 &   5.07 &   \color{blue}{\textbf{0.25}} &   \color{blue}{0.26} & \color{blue}{0.34} &   \color{blue}{0.33} &  228.1 &   127.1 &   1.96 &  1.99 &   1.95 &  \textbf{1.88}  \\ \hline
       & 20 dB & 8.11 &   8.10 &   4.50 & \textbf{3.98} & 6.58 &   6.16 &  174.6 &   165.1 &   121.2 &   122.0 &   \textbf{120.8} & 121.8 \\
      OCO-2 & 40 dB &  8.04 &   7.80 &   \color{blue}{0.79} &   \color{blue}{\textbf{0.73}} &   1.12 &   1.04 & 46.35 &  38.60 &   \textbf{2.33} &   2.35 &   2.34 &   2.35 \\
      Band 1 & 55 dB &  8.03 &   7.79 &  \color{blue}{\textbf{0.54}} & \color{blue}{0.56} & \color{blue}{0.84} &   \color{blue}{0.76} &   45.11 &  37.36 &   \textbf{1.15} &   1.17 &   1.18 & 1.16 \\   
       & 80 dB &  8.03 &   7.79 &   \color{blue}{\textbf{0.52}} & \color{blue}{0.55} & \color{blue}{0.83} & \color{blue}{0.73} &  45.07 &  37.32 & \textbf{1.11} &   1.13 &   1.14 & 1.13
  \\ 
       & 120 dB & 8.03 &   7.79 &   \color{blue}{\textbf{0.51}} & \color{blue}{0.55} & \color{blue}{0.83} &   \color{blue}{0.73} &  45.07 & 37.32 &   \textbf{1.11} &  1.13 &   1.14 &   1.13 \\ \hline
       & 20 dB & 16.28 &  \textbf{3.39} &   4.58 & 4.37 & 14.38 &  14.23 &  185.5 &  116.4 &   \textbf{112.3} &   112.4 &   112.7 & 112.9
  \\
       MicroCarb & 40 dB & 16.27 &   2.04 &   \color{blue}{\textbf{0.54}} & \color{blue}{0.56} & 2.05 &   2.37 & 70.2 &   1.56 &   \textbf{1.32} &   \textbf{1.32} &   1.53 & 1.70  \\
      Band 1 & 55 dB & 16.27 &   2.03 &   \color{blue}{\textbf{0.29}} & \color{blue}{0.33} &   1.33 &   1.66 &  69.21 &   0.43 &   \textbf{0.23} &   0.24 & 0.38 &   0.61 \\   
       & 80 dB & 16.27 &   2.03 & \color{blue}{\textbf{0.28}} &   \color{blue}{0.32} & 1.27 &  1.68 &  69.2 &   0.39 &   \textbf{0.20}&  \textbf{0.20} &   0.34 & 0.60 \\ 
       & 120 dB & 16.27 &  2.03 & \color{blue}{\textbf{0.28}} &   \color{blue}{0.32} & 1.27 &  1.69 & 69.2 &   0.39 &   \textbf{0.20} &   \textbf{0.20} &   0.34 & 0.60 \\\hline
    \end{tabular}
    \end{center}

\end{table}

\subsubsection{Parameter tuning for SPIRIT}
\begin{figure}[h!] 
 \includegraphics[trim={4.5cm 0cm 0cm 0cm}, clip,width=18.7cm]{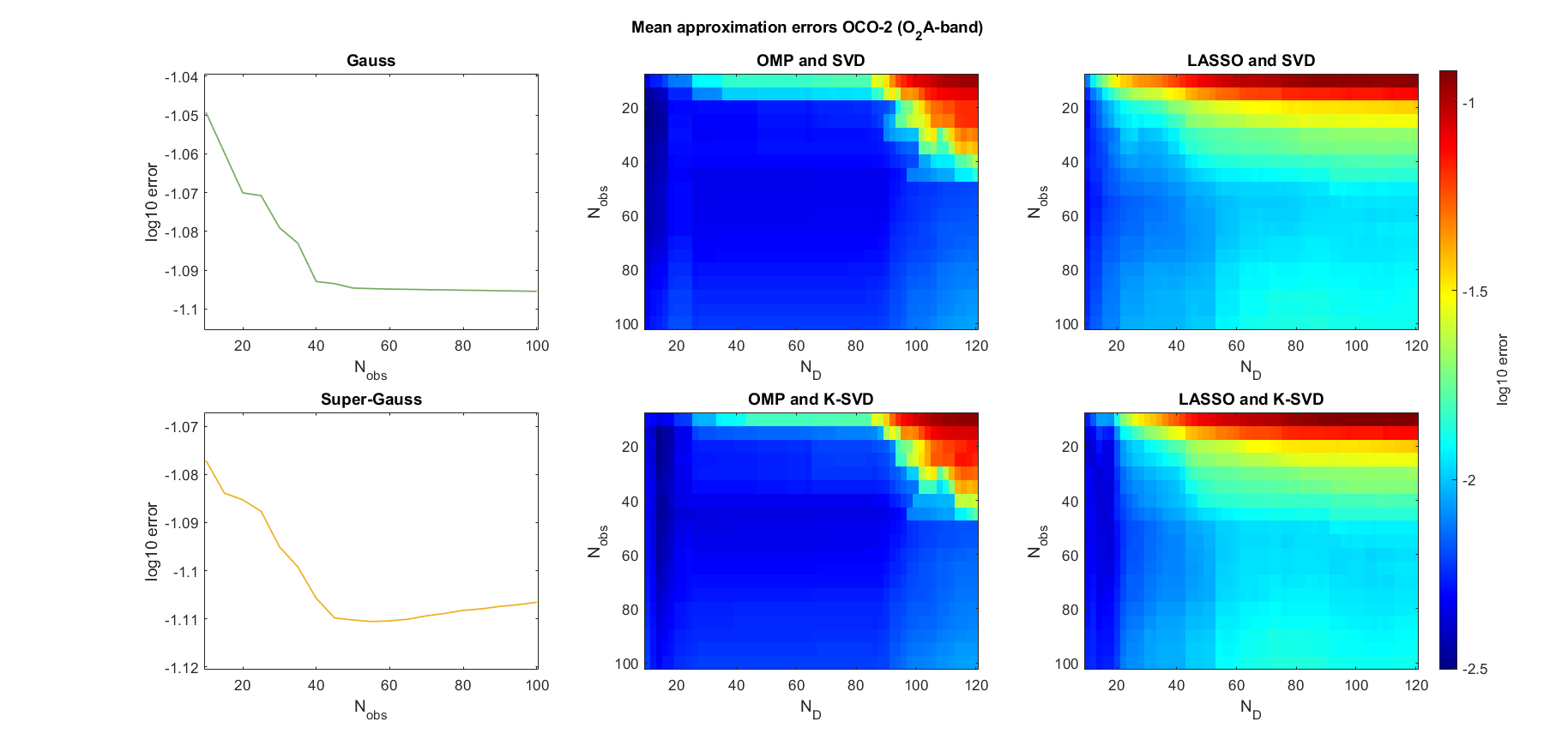}
 \caption{Mean approximation errors for OCO-2 and the different estimation methods (Gauss, Super-Gauss, OMP and LASSO with SVD or K-SVD) versus the number of observations $N_{\text{obs}}$ and the dictionary size $N_{\text{D}}$ for K = 5.} \label{OCO}
\end{figure}

\begin{figure}[h!] 
 \includegraphics[trim={4.2cm 0cm 0cm 0cm}, clip,width=18.7cm]{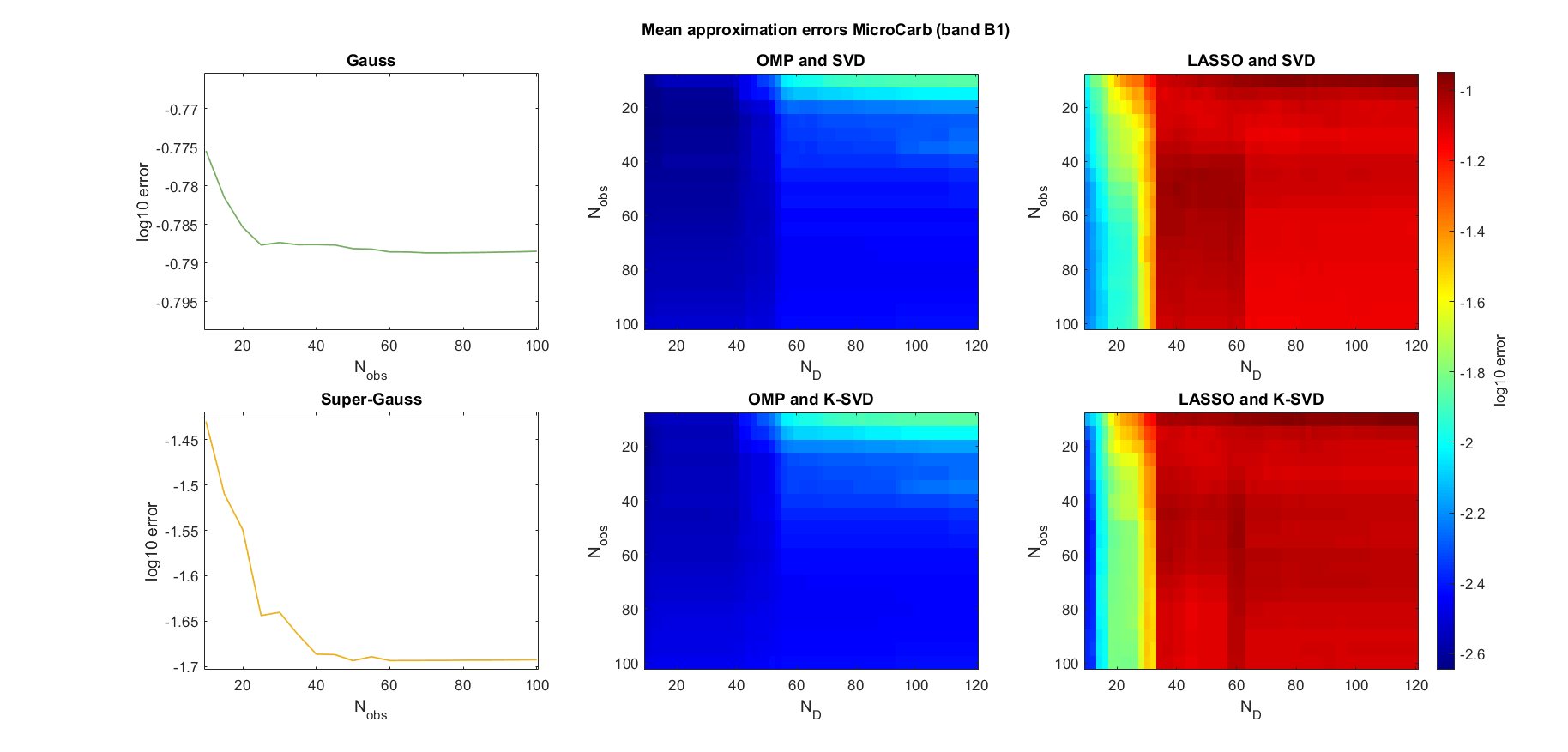}
 \caption{Mean approximation errors for MicroCarb and the different estimation methods (Gauss, Super-Gauss, OMP and LASSO with SVD or K-SVD) versus the number of observations $N_{\text{obs}}$ and the dictionary size $N_{\text{D}}$ for K = 4.} \label{MC}
\end{figure}

The proposed approach requires the choice of a small number of parameters, namely the size of the sliding window $N_{\text{obs}}$, the size of the dictionary $N_{\text{D}}$ and the number of atoms $K$.
The number of atoms $K$ has been discussed above and found to yield robust results for $K\approx 3-5$, across all instruments, see Figs.~\ref{figure:Avantes_GOME-2_spectrometer} to \ref{figure:OCO_2_CNES_spectrometer} and the corresponding discussions in Section \ref{sec:perfwav}.
This section further studies the impact of $N_{\text{obs}}$ and $N_{\text{D}}$ on the ISRF approximation errors. To this end, Figs.~\ref{OCO} and \ref{MC} show the approximation errors (in log10 scale) as a function of $N_{\text{obs}}$ for the Gaussian and Super-Gaussian parameterizations, and as functions of $(N_{\text{obs}}, N_{\text{D}})$  for the methods based on sparse representations. This analysis is conducted for the two instruments OCO-2 and MicroCarb without loss of generality and the results are averaged for all ISRFs. 

The ISRF estimation errors decrease as $N_{\text{obs}}$ increases, as expected. However, this decrease is more important for sparse methods (e.g., for $N_{\text{obs}}=80$, the mean ISRF errors for Gauss, Super-Gauss and SPIRT are equal to 16.27\%, 2.04\%, whereas they are equal to 0.29\% for OMP/SVD, 0.33\% for OMP/K-SVD, 1.23\% for LASSO/SVD and 1.40\% for LASSO/K-SVD, showing the interest of exploiting sparsity for ISRF estimation. Note that the results displayed in Figs. \ref{OCO} and \ref{MC} show that it is not recommended to use too many atoms in the dictionaries since the ISRF estimation errors are relatively high for large values of $N_{\text{D}}$. Based on this observation, the value $N_{\text{D}}=25$ was used in all the other experiments.

\subsubsection{Robustness to ISRF changes}
\label{Subsubsec:robustnesstoISRFchanges}
\begin{figure}[h!] 
\centering 
            \includegraphics[trim={0cm 0cm 0cm 0cm}, clip,width=10cm]{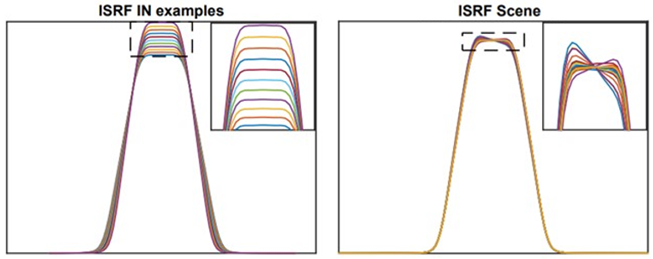}
            \caption{Examples of ISRFs IN, scene and NU (MicroCarb band 1).}
    \label{figure:ISRF_change}
\end{figure}

The ISRFs considered in the previous sections were obtained from uniform scenes referred to as “ISRF IN” for the MicroCarb mission (see Fig.~\ref{figure:ISRF_change} for illustrations). However these ISRFs can change depending on the scene observed by the instrument.

\begin{figure}[h!]
    \begin{minipage}{0.01\textwidth}
       \centering
       \vspace*{-1em}
       \hspace*{1cm}
            \includegraphics[trim={0cm 0cm 0cm 0cm}, clip,width=8cm]{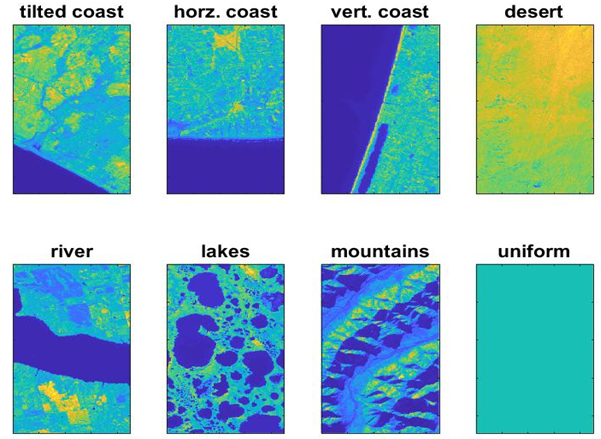}
    \end{minipage}\hfill\begin{minipage}{0.45\textwidth}
        \centering
            \includegraphics[trim={0cm 0cm 0cm 0cm}, clip,width=7cm]{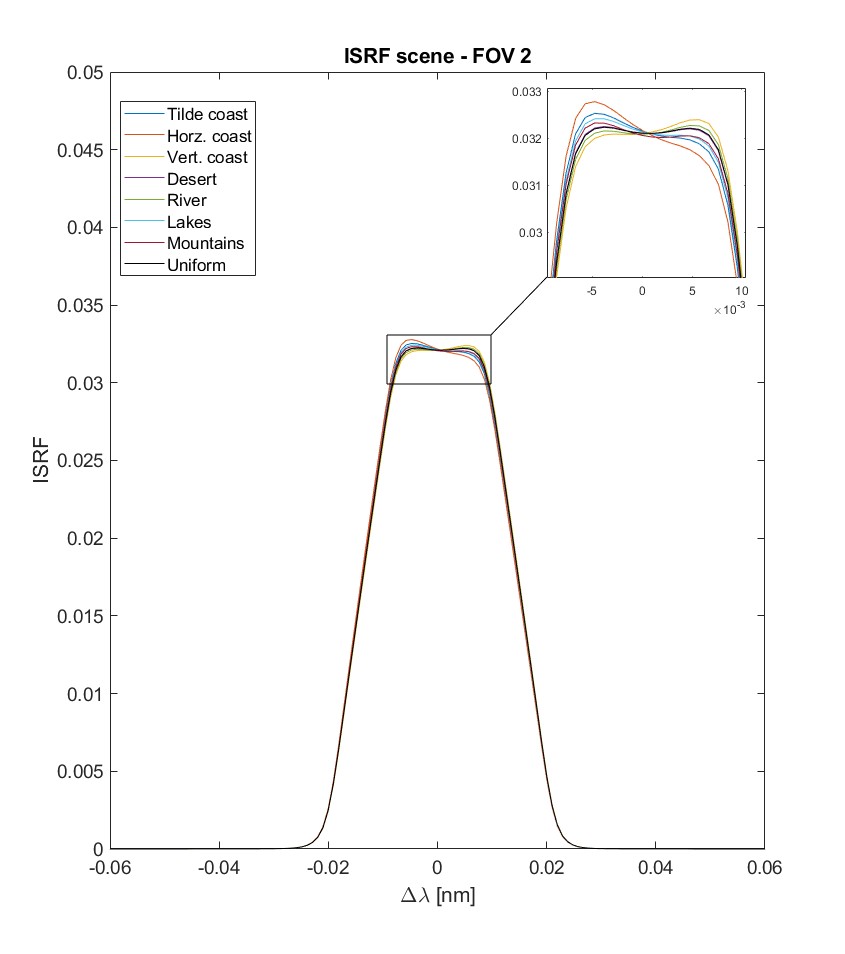}
    \end{minipage} 
     \caption{Eight types of scenes (left) with the corresponding ISRFs (FOV 2) (right) for the MicroCarb instrument.}
    \label{figure:ISRF_change_scenes}
\end{figure}

\begin{figure}[h!] 
\centering   
            \includegraphics[trim={1cm 0.5cm 1cm 0.5cm}, clip,width=19cm]{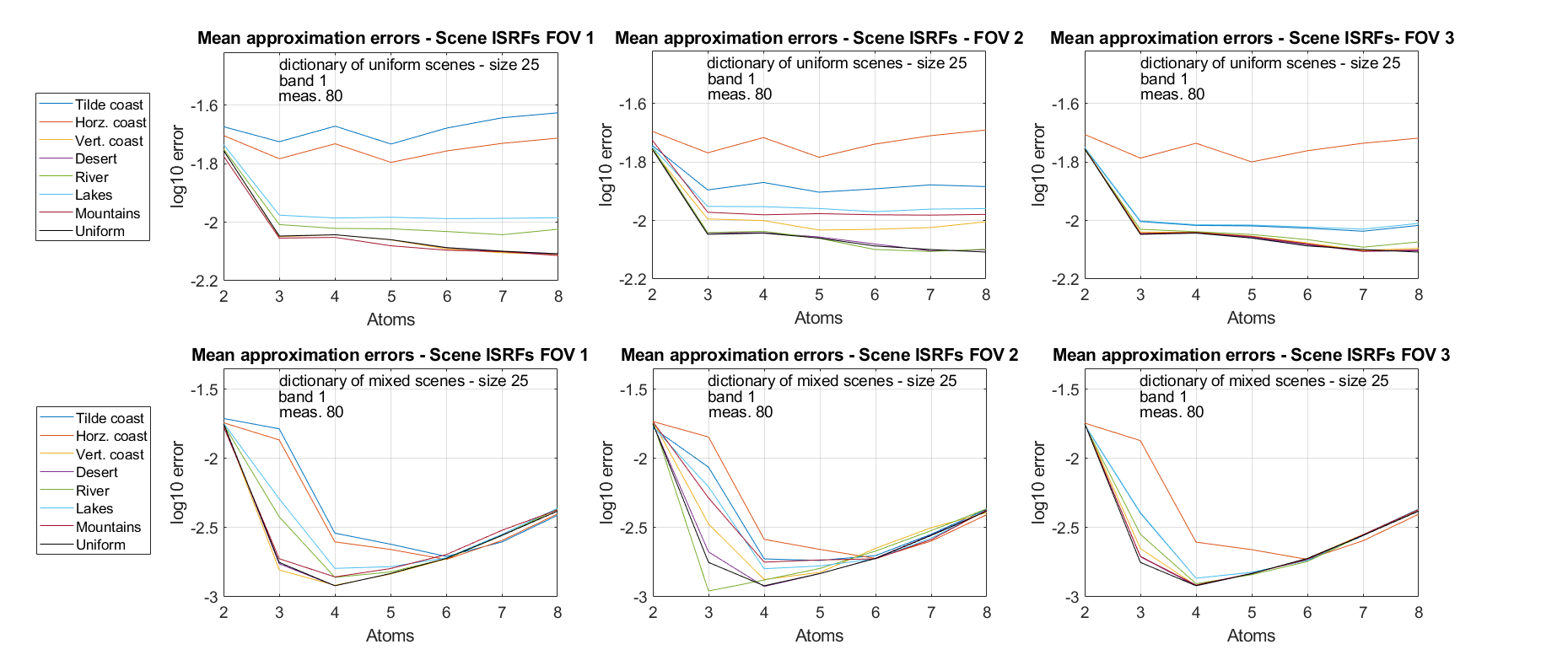}
            \caption{ISRF estimation errors for Scene ISRFs obtained using a dictionary of uniform ISRFs (top) and of mixed ISRFs (bottom).}
    \label{figure:results_ISRF_change_mixedDico}
\end{figure}

The design of the MicroCarb instrument makes the ISRF sensitive to the slit illumination during the integration time. If the slit is not uniformly illuminated, the ISRF shape changes as shown in Fig.~\ref{figure:ISRF_change_scenes}, which displays an example of eight different scenes of the Earth’s surface with their associated ISRFs. Note that in this figure, the ISRF of a desert scene is very similar to the ISRF of a uniform scene, contrary to the ISRF of an horizontal coast profil, which makes the slit blinded during a third of the integration time and leads to an asymmetric left-distorted ISRF. For each of these eight scenes, three Fields Of Views (FOVs) were defined and labeled as FOV1, FOV2 and FOV3 in Fig.~\ref{figure:results_ISRF_change_mixedDico}). These FOVs correspond to a binning into three parts of the detector pixels along the ACT axis. The spatial pixels in each FOV are averaged to increase the spectral SNR. This binning and averaging step allow three measured spectra per imaged area to be determined, whose ISRFs have to be estimated.

This section first studies the performance of SPIRIT for estimating non uniform Scene ISRFs for the first band (band B1) of the MicroCarb spectrometer using a dictionary learned from examples of uniform ISRFs (ISRF IN). The top part of Fig.~\ref{figure:results_ISRF_change_mixedDico} presents the results obtained by SPIRIT (using the OMP algorithm and SVD of the on-ground ISRFs). Estimating the ISRFs Scene with basis functions calculated from uniform ISRFs is clearly a more difficult problem. Indeed, the resulting normalized ISRF errors exceed 1 \% for several scenes and FOVs. In order to take into account the diversity of ISRFs, a new dictionary constructed from several examples of ISRF IN and ISRF Scene can be considered. More precisely, a matrix composed of $103$ ISRFs IN (one out of ten) and $3$ ISRFs Scene (out of $24$), i.e., of size $(103+3 = 106) \times N$ (with $N=895$) was considered to build of dictionary of $N_D = 25$ atoms. The bottom part of Fig. \ref{figure:results_ISRF_change_mixedDico} shows the performance obtained for the first MicroCarb band for Scene ISRFs. Using only three additional examples of ISRFs Scene in the dictionary construction step allows ISRF estimation errors to be less than 1\% criterion, which is the target for the estimation accuracy. Note that the lowest approximation errors are obtained in most cases using three to $K=6$ atoms from the dictionary, as before. To conclude, this section showed that including diversity in the ISRFs used for the construction of the dictionary, leads to significant gains in ISRF estimation.

\clearpage

\conclusions 
This article proposed a new method for estimating the instrument spectral response functions (ISRFs) of spectrometers. This method was based on a sparse decomposition of these spectral responses into a dictionary of basis functions called atoms. The resulting new estimation strategy was applied to a large variety of ISRFs associated with different spectrometers such as Avantes, GOME-2, OMI, TROPOMI, OCO-2 and Microcarb. The response functions of these instruments can be estimated with high accuracy with the proposed method, leading to normalized errors of less than 1\%. The proposed estimation method decomposed the ISRFs of interest on a dictionary of reference atoms constructed using an SVD decomposition of ISRF examples or the K-SVD algorithm.
This decomposition was conducted using two different sparse coding algorithms based on orthogonal matching pursuit (OMP) and the least absolute shrinkage and selection operator (LASSO). Overall, the results obtained using the OMP algorithm and the SVD are very promising, yielding small estimation errors and fast execution times, which suggests its practicality for in-flight scenarios. Another interesting property of the proposed estimation method is that is not impacted significantly by the shape of the ISRF, allowing accurate estimations for different types of scenes.

Future work includes the consideration of radiometric and spectral errors (such as straylight, residual errors of calibration or spectral shifts) that can degrade the performance of ISRF estimation. These errors are expected to affect more significantly some specific wavelengths, which might deserve to investigate specific algorithms jointly correcting the errors and estimating the ISRFs. The resulting problem is more challenging since there are non-linear relationships between the spectrometer measurements and these radiometric and spectral errors. Another interesting prospect is to study the potential interest of machine learning algorithms for error correction and ISRF estimation.

\codedataavailability{More details on the data and code used in this study are available upon request from the corresponding author.}

\noappendix      

\appendixfigures 

\appendixtables

\authorcontribution{JE gathered the data for different spectrometers and CP for the MicroCarb spectrometer. JMG and CP contributed to a first formalization of the problem. Mathematical formulation, implementation and formal analysis were conducted by JE, JYT and HW.  All authors have contributed to the writing process through discussion and feedback.} 

\competinginterests{The authors declare that they have no conflict of interest.}

\begin{acknowledgements}
This study was supported by the French Space Agency (CNES), France and by Thales Alenia Space Cannes, France. We would like to thank M. Denis Jouglet at the department of Atmospheric Sounding at CNES for providing the reference spectra at the different wavelengths. Moreover, we express our gratitude to Steffen Beirle at the Max Planck Institute for Chemistry(MPI-C) for helpful discussions and for provided some ISRF data. 
\end{acknowledgements}


\bibliographystyle{copernicus}
\bibliography{biblio.bib}

\begin{thebibliography}{29}
\providecommand{\natexlab}[1]{#1}
\providecommand{\url}[1]{\texttt{#1}}
\providecommand{\urlprefix}{}
\expandafter\ifx\csname urlstyle\endcsname\relax
  \providecommand{\doi}[1]{https://doi.org/\discretionary{}{}{}#1}\else
  \providecommand{\doi}{https://doi.org/\discretionary{}{}{}\begingroup
  \urlstyle{rm}\Url}\fi

\bibitem[{Ava()}]{Avantes}
{AvaSpec-ULS2048x64-EVO} - {Data User's product description},
  \urlprefix\url{https://avantes.b-cdn.net/content/uploads/2020/11/DS-Spec-Avaspec-ULS2048x64-EVO-200702.pdf
  (last access: 06 November 2023)}.

\bibitem[{Aharon et~al.(2006)Aharon, Elad, and Bruckstein}]{Elad2006}
Aharon, M., Elad, M., and Bruckstein, A.: {K-SVD}: An algorithm for designing
  overcomplete dictionaries for sparse representation, IEEE Trans. Signal
  Process., 54, 4311--4322, \doi{10.1109/TSP.2006.881199}, 2006.

\bibitem[{Armante et~al.(2013)Armante, Scott, Capelle, Chédin, Bernard,
  Standfuss, Tournier, and Pierangelo}]{4AOPPresentation}
Armante, R., Scott, N.~A., Capelle, V., Chédin, A., Bernard, E., Standfuss,
  C., Tournier, B., and Pierangelo, C.: {IASI conference Presentation : 4AOP :
  A fast and accurate operational forward radiative transfer model},
  \urlprefix\url{https://cnes.fr/sites/default/files/migration/smsc/iasi/PDF/conf3/posters/84_Armante_R.pdf},
  2013.

\bibitem[{Babic et~al.(2022)Babic, Braak, Dierssen, Kissi-Ameyaw, Kleipool,
  Leloux, Loots, Ludewig, Rozemeijer, Smeets, and Vacanti}]{Tropomi_document}
Babic, L., Braak, R., Dierssen, W., Kissi-Ameyaw, J., Kleipool, Q., Leloux, J.,
  Loots, E., Ludewig, A., Rozemeijer, N., Smeets, J., and Vacanti, G.:
  {Algorithm theoretical basis document for the TROPOMI L01b data processor},
  \urlprefix\url{sentinels.copernicus.eu/documents/247904/2476257/Sentinel-5P-TROPOMI-Level-1B-ATBD},
  2022.

\bibitem[{Beirle et~al.(2017)Beirle, Lampel, Lerot, Sihler, and
  Wagner}]{Beirle2017}
Beirle, S., Lampel, J., Lerot, C., Sihler, H., and Wagner, T.: {Parameterizing
  the instrumental spectral response function and its changes by a
  super-Gaussian and its derivatives}, Atmos. Meas. Tech., 10, 581--598,
  \doi{10.5194/amt-10-581-2017}, 2017.

\bibitem[{Bertaux et~al.(2020)Bertaux, Hauchecorne, Lef\`evre, Br\'eon, Blanot,
  Jouglet, Lafrique, and Akaev}]{Bertaux2020}
Bertaux, J.-L., Hauchecorne, A., Lef\`evre, F., Br\'eon, F.-M., Blanot, L.,
  Jouglet, D., Lafrique, P., and Akaev, P.: The use of the 1.27\,\unit{\mu}m
  \chem{O_2} absorption band for greenhouse gas monitoring from space and
  application to MicroCarb, Atmospheric Measurement Techniques, 13, 3329--3374,
  \doi{10.5194/amt-13-3329-2020}, 2020.

\bibitem[{Cansot et~al.(2022)Cansot, Pistre, Castelnau, Landiech, Georges,
  Gaeremynck, and Bernard}]{Cansot2022}
Cansot, E., Pistre, L., Castelnau, M., Landiech, P., Georges, L., Gaeremynck,
  Y., and Bernard, P.: MICROCARB INSTRUMENT, OVERVIEW AND FIRST RESULTS, Proc.
  SPIE 12777, Inf. Conf. Space Optics, 12777, 1--13, \doi{10.1117/12.2690330},
  2022.

\bibitem[{Crisp et~al.(2021)Crisp, Rosenberg, Chapsky, Keller~Rodrigues, Lee,
  Merrelli, Osterman, Oyafuso, Pollock, Spiers, Yu, Zong, and
  Eldering}]{OCO2_document}
Crisp, D., Rosenberg, R., Chapsky, L., Keller~Rodrigues, G.~R., Lee, R.,
  Merrelli, A., Osterman, G., Oyafuso, F., Pollock, R., Spiers, G., Yu, S.,
  Zong, J., and Eldering, A.: {Orbiting Carbon Observatory – 2 \& 3 (OCO-2 \&
  OCO-3)- Level 1B! Algorithm Theoretical Basis Document},
  \urlprefix\url{sentinels.copernicus.eu/documents/247904/2476257/Sentinel-5P-TROPOMI-Level-1B-ATBD},
  2021.

\bibitem[{Dobber et~al.(2006)Dobber, Dirksen, Levelt, van~den Oord, Voors,
  Kleipool, Jaross, Kowalewski, Hilsenrath, Leppelmeier, de~Vries, Dierssen,
  and Rozemeijer}]{Dobber2006}
Dobber, M.~R., Dirksen, R.~J., Levelt, P.~F., van~den Oord, G. H.~J., Voors, R.
  H.~M., Kleipool, Q., Jaross, G., Kowalewski, M., Hilsenrath, E., Leppelmeier,
  G.~W., de~Vries, J., Dierssen, W., and Rozemeijer, N.~C.: Ozone monitoring
  instrument calibration, IEEE Trans. Geosci. Remote Sens., 44, 1209--1238,
  \doi{10.1109/TGRS.2006.869987}, 2006.

\bibitem[{Figueiredo et~al.(2007)Figueiredo, Nowak, and
  Wright}]{Gradient_Projection}
Figueiredo, M. A.~T., Nowak, R.~D., and Wright, S.~J.: Gradient Projection for
  Sparse Reconstruction: Application to Compressed Sensing and Other Inverse
  Problems, IEEE J. Sel. Topics Signal Process., 1, 586--597,
  \doi{10.1109/JSTSP.2007.910281}, 2007.

\bibitem[{Kim et~al.(2007)Kim, Koh, Lustig, Boyd, and
  Gorinevsky}]{Interior-Point}
Kim, S.-J., Koh, K., Lustig, M., Boyd, S., and Gorinevsky, D.: An
  Interior-Point Method for Large-Scale $\ell_1$-Regularized Least Squares,
  IEEE J. Sel. Topics Signal Process., 1, 606--617,
  \doi{10.1109/JSTSP.2007.910971}, 2007.

\bibitem[{Kleipool et~al.(2018)Kleipool, Ludewig, Babi\'c, Bartstra, Braak,
  Dierssen, Dewitte, Kenter, Landzaat, Leloux, Loots, Meijering, van~der Plas,
  Rozemeijer, Schepers, Schiavini, Smeets, Vacanti, Vonk, and
  Veefkind}]{Kleipool2018}
Kleipool, Q., Ludewig, A., Babi\'c, L., Bartstra, R., Braak, R., Dierssen, W.,
  Dewitte, P.-J., Kenter, P., Landzaat, R., Leloux, J., Loots, E., Meijering,
  P., van~der Plas, E., Rozemeijer, N., Schepers, D., Schiavini, D., Smeets,
  J., Vacanti, G., Vonk, F., and Veefkind, P.: Pre-launch calibration results
  of the TROPOMI payload on-board the Sentinel-5 Precursor satellite, Atmos.
  Meas. Tech., 11, 6439--6479, \doi{10.5194/amt-11-6439-2018}, 2018.

\bibitem[{Lagarias et~al.(1998)Lagarias, Reeds, Wright, and
  Wright}]{Nelder_Mead}
Lagarias, J.~C., Reeds, J.~A., Wright, M.~H., and Wright, P.~E.: Convergence
  Properties of the {N}elder-{M}ead Simplex Method in Low Dimensions, SIAM J.
  Optim., 9, 112--147, 1998.

\bibitem[{Lee et~al.(2017)Lee, O’Dell, Wunch, Roehl, Osterman, Blavier,
  Rosenberg, Chapsky, Frankenberg, Hunyadi-Lay, Fisher, Rider, Crisp, and
  Pollock}]{Lee2017}
Lee, R. A.~M., O’Dell, C.~W., Wunch, D., Roehl, C.~M., Osterman, G.~B.,
  Blavier, J.-F., Rosenberg, R., Chapsky, L., Frankenberg, C., Hunyadi-Lay,
  S.~L., Fisher, B.~M., Rider, D.~M., Crisp, D., and Pollock, R.: Preflight
  Spectral Calibration of the Orbiting Carbon Observatory 2, IEEE Trans.
  Geosci. Remote Sens., 55, 2499--2508, \doi{10.1109/TGRS.2016.2645614}, 2017.

\bibitem[{Mallat and Zhang(1993)}]{Mallat1993}
Mallat, S.~G. and Zhang, Z.: Matching pursuits with time-frequency
  dictionaries, IEEE Trans. Signal Process., 41, 3397--3415,
  \doi{10.1109/78.258082}, 1993.

\bibitem[{Munro et~al.(2016)Munro, Lang, Klaes, Poli, Retscher, Lindstrot,
  Huckle, Lacan, Grzegorski, Holdak, Kokhanovsky, Livschitz, and
  Eisinger}]{Munro2016}
Munro, R., Lang, R., Klaes, D., Poli, G., Retscher, C., Lindstrot, R., Huckle,
  R., Lacan, A., Grzegorski, M., Holdak, A., Kokhanovsky, A., Livschitz, J.,
  and Eisinger, M.: {The GOME-2 instrument on the Metop series of satellites:
  instrument design, calibration, and level 1 data processing -- an overview},
  Atmos. Meas. Tech., 9, 1279--1301, \doi{10.5194/amt-9-1279-2016}, 2016.

\bibitem[{{NOVELTIS} et~al.(){NOVELTIS}, {CNES}, and {LMD}}]{4AOP}
{NOVELTIS}, {CNES}, and {LMD}: {4A/OP} - Operational release for {4A} -
  {Automatized Atmospheric Absorption Atlas},
  \urlprefix\url{https://4aop.noveltis.fr/references-and-publications}.

\bibitem[{Pati et~al.(1993)Pati, Rezaiifar, and Krishnaprasad}]{Pati1993}
Pati, Y.~C., Rezaiifar, R., and Krishnaprasad, P.~S.: Orthogonal matching
  pursuit: recursive function approximation with applications to wavelet
  decomposition, in: Proc. Asilomar Conf. Signals, Systems and Computers, pp.
  40--44, Pacific Grove, CA, USA, \doi{10.1109/ACSSC.1993.342465}, 1993.

\bibitem[{Schenkeveld et~al.(2017)Schenkeveld, Jaross, Marchenko, Haffner,
  Kleipool, Rozemeijer, Veefkind, and Levelt}]{Schenkeveld2017}
Schenkeveld, V. M.~E., Jaross, G., Marchenko, S., Haffner, D., Kleipool, Q.~L.,
  Rozemeijer, N.~C., Veefkind, J.~P., and Levelt, P.~F.: {In-flight performance
  of the Ozone Monitoring Instrument}, Atmos. Meas. Tech., 10, 1957--1986,
  \doi{10.5194/amt-10-1957-2017}, 2017.

\bibitem[{Siddans and Latter(2018)}]{Siddans2018}
Siddans, R. and Latter, B.~G.: {Analysis of GOME-2 (FM201-3): Slit function
  Measurements Final Report Eumetsat Contract No. EUM/CO/04/1298/RM},
  \urlprefix\url{https://www.eumetsat.int/media/42779}, 2018.

\bibitem[{Smeets et~al.(2002)Smeets, Kleipool, van Hees, and
  M.}]{Tropomi_dataUser}
Smeets, J., Kleipool, Q., van Hees, R., and M., S.: {Readme for TROPOMI
  instrument spectral response functions},
  \urlprefix\url{https://sentinels.copernicus.eu/documents/247904/3541451/S5P-KNMI-OCAL-0149-ME-readme_for_TROPOMI_UVN_instrument_spectral_response_function-3.0.0-20180401_0.pdf},
  2002.

\bibitem[{Sun et~al.(2017{\natexlab{a}})Sun, Liu, Huang, Gonz\'alez~Abad, Cai,
  Chance, and Yang}]{Sun2017a}
Sun, K., Liu, X., Huang, G., Gonz\'alez~Abad, G., Cai, Z., Chance, K., and
  Yang, K.: Deriving the slit functions from {OMI} solar observations and its
  implications for ozone-profile retrieval, Atmos. Meas. Tech., 10, 3677--3695,
  \doi{10.5194/amt-10-3677-2017}, 2017{\natexlab{a}}.

\bibitem[{Sun et~al.(2017{\natexlab{b}})Sun, Liu, Nowlan, Cai, Chance,
  Frankenberg, Lee, Pollock, Rosenberg, and Crisp}]{Sun2017b}
Sun, K., Liu, X., Nowlan, C.~R., Cai, Z., Chance, K., Frankenberg, C., Lee, R.
  A.~M., Pollock, R., Rosenberg, R., and Crisp, D.: Characterization of the
  OCO-2 instrument line shape functions using on-orbit solar measurements,
  Atmos. Meas. Tech., 10, 939--953, \doi{10.5194/amt-10-939-2017},
  2017{\natexlab{b}}.

\bibitem[{Tan et~al.(2015)Tan, Tsang, and Wang}]{LASSO_PART1}
Tan, M., Tsang, I.~W., and Wang, L.: Matching Pursuit {LASSO Part I}: Sparse
  Recovery Over Big Dictionary, IEEE Trans. Signal Process., 63, 727--741,
  \doi{10.1109/TSP.2014.2385036}, 2015.

\bibitem[{Tibshirani(1996)}]{tibshirani96regression}
Tibshirani, R.: Regression Shrinkage and Selection via the Lasso, Journal of
  the Royal Statistical Society (Series B), 58, 267--288, 1996.

\bibitem[{Tošić and Frossard(2011)}]{dictionary_learning}
Tošić, I. and Frossard, P.: Dictionary Learning, IEEE Signal Process. Mag.,
  28, 27--38, \doi{10.1109/MSP.2010.939537}, 2011.

\bibitem[{van Hees et~al.(2018)van Hees, Tol, Cadot, Krijger, Persijn, van
  Kempen, Snel, Aben, and Hoogeveen}]{vanHees2018}
van Hees, R.~M., Tol, P. J.~J., Cadot, S., Krijger, M., Persijn, S.~T., van
  Kempen, T.~A., Snel, R., Aben, I., and Hoogeveen: Determination of the
  TROPOMI-SWIR instrument spectral response function, Atmospheric Measurement
  Techniques, 11, 3917--3933, \doi{10.5194/amt-11-3917-2018}, 2018.

\bibitem[{Zhang and Huang(2008)}]{LASSO_LS}
Zhang, C.-H. and Huang, J.: {The sparsity and bias of the Lasso selection in
  high-dimensional linear regression}, Ann. Stat., 36, 1567 -- 1594,
  \doi{10.1214/07-AOS520}, 2008.

\bibitem[{Zhang et~al.(2015)Zhang, Xu, Yang, Li, and Zhang}]{SurveySparseRep}
Zhang, Z., Xu, Y., Yang, J., Li, X., and Zhang, D.: A Survey of Sparse
  Representation: Algorithms and Applications, IEEE Access, 3, 490--530,
  \doi{10.1109/ACCESS.2015.2430359}, 2015.

\end{thebibliography}

\end{document}